\documentclass[10pt]{amsart}
\usepackage{amssymb, amsthm, epsfig}
\usepackage[mathscr]{euscript}
\usepackage{pb-diagram}

\newtheorem{thm}{Theorem}[section]
\newtheorem{cor}[thm]{Corollary}
\newtheorem{lem}[thm]{Lemma}
\newtheorem{prop}[thm]{Proposition}
\newtheorem{conj}[thm]{Conjecture}

\newtheorem{defnn}[thm]{Definition}
\newenvironment{definition}{\begin{defnn} \em}{\end{defnn}}
\newtheorem{remarkk}[thm]{Remark}
\newenvironment{remark}{\begin{remarkk} \em}{\end{remarkk}}
\newtheorem{examplee}[thm]{Example}

\newcommand{\scrc}{\mathcal{C}}
\newcommand{\scrdc}{\mathcal{D}\mathcal{C}}
\newcommand{\scrf}{\mathcal{F}}
\newcommand{\g}{\mathcal{G}}
\newcommand{\scrk}{\mathcal{K}}

\newcommand{\bbz}{\mathbb{Z}}

\newcommand{\bbq}{\mathbb{Q}}
\newcommand{\bbn}{\mathbb{N}}

\newcommand{\G}{\Gamma}

\newcommand{\Hom}{\operatorname{Hom}}
\newcommand{\Ker}{\operatorname{Ker}}
\def\L{\Lambda}
\def\G{\Gamma}
\def\e{\epsilon}
\def\int{\mathrm{int}}

\def\Im{\mathrm{Im}}

\def\Hom{\mathrm{Hom}}
\def\S{\Sigma}

\title{New obstructions to doubly slicing knots}
\author{Taehee Kim}
\date{\today}
\address{Math Department--MS 136, Rice University, 6100 S.~Main
St.~, Houston, TX 77005-1892, USA} \email{tkim@rice.edu}
\def\subjclassname{\textup{2000} Mathematics Subject Classification}
\expandafter\let\csname
subjclassname@1991\endcsname=\subjclassname
\expandafter\let\csname
subjclassname@2000\endcsname=\subjclassname 

\begin{document}
\begin{abstract}
A knot in the 3-sphere is called doubly slice if it is a slice of
an unknotted 2-sphere in the 4-sphere. We give a bi-sequence of
new obstructions for a knot being doubly slice. We construct it
following the idea of Cochran-Orr-Teichner's filtration of the
classical knot concordance group. This yields a bi-filtration of
the monoid of knots (under the connected sum operation) indexed by
pairs of half integers. Doubly slice knots lie in the intersection
of this bi-filtration. We construct examples of knots which
illustrate non-triviality of this bi-filtration at all levels. In
particular, these are new examples of algebraically doubly slice
knots that are not doubly slice, and many of these knots are
slice. Cheeger-Gromov's von Neumann rho invariants play a key role
to show non-triviality of this bi-filtration. We also show that
some classical invariants are reflected at the initial levels of
this bi-filtration, and obtain a bi-filtration of the double
concordance group.
\end{abstract}

\maketitle \vspace*{-2em}

\section{Introduction}
\label{sec:introduction} We work in the topologically locally flat
category. An $n$-knot in the $(n+2)$-sphere is called {\em doubly
slice} (or {\em doubly null cobordant}) if it is a slice of an
unknotted $(n+1)$-sphere in the $(n+3)$-sphere. The notion of
doubly slice knots was introduced by Fox \cite{Fo} in the 60's.
For odd dimensional knots, Sumners \cite{Su} showed that if a knot
is doubly slice, then it has an associated Seifert form which is
hyperbolic. We call the knots satisfying this Seifert form
condition {\em algebraically doubly slice} (or {\em algebraically
doubly null cobordant}). It was shown that for odd high
dimensional simple knots, this Seifert form obstruction is
sufficient for being doubly slice \cite{Su, Ke2}. This result was
generalized to even high dimensional knots by Stoltzfus \cite{S1,
S2} using the obstructions based on the linking form defined by
Levine \cite{L} and Farber \cite{F}. In this paper, we work in the
classical dimension. So by ``knot" we mean a 1-knot in the
3-sphere unless mentioned otherwise.

In \cite{GL}, Gilmer and Livingston showed that there exists a
slice knot which is algebraically doubly slice but not doubly
slice. (A knot is called {\em slice} if it bounds a locally flat
2-disk in the 4-ball.) One can see that if a knot is doubly slice
then every finite branched cyclic cover of the knot is embedded in
the 4-sphere. They applied their own obstructions to embedding
3-manifolds into the 4-sphere to show that their example is not
doubly slice. High dimensional analogues of this result were
obtained by Ruberman \cite{Ru}. Recently Friedl \cite{Fri} found
doubly slicing obstructions using eta invariants associated to
finite dimensional unitary representations.

Meanwhile, Cochran, Orr, and Teichner (henceforth COT) established
a filtration of the classical knot concordance group $\scrc$
\cite{COT1}.
\[
 0\subset\cdots\subset\scrf_{n.5}\subset\scrf_{n}\subset\cdots
 \subset\scrf_{1.5}\subset\scrf_{1.0}\subset\scrf_{0.5}
 \subset\scrf_{0}\subset\scrc
\]
where $\scrf_{m}$ is the set of all $(m)$-solvable knots. Roughly
speaking, a 3-manifold is said to be {\em $(m)$-solvable} (via
$W$) if it bounds a spin 4-manifold $W$ that induces an
isomorphism on the first homology and satisfies a certain
condition on the intersection form of the $m^{\text{th}}$ derived
cover of $W$. A knot is called \emph{$(m)$-solvable} (via $W$) if
zero surgery on the knot in the 3-sphere is $(m)$-solvable (via
$W$). If $K$ is $(m)$-solvable via $W$, then $W$ is called an {\em
$(m)$-solution} for the knot (or for zero surgery on the knot in
the 3-sphere). COT showed that if a knot is $(1.5)$-solvable, then
all the previously known concordance invariants including
Casson-Gordon invariants vanish for the knot \cite[Theorem
9.11]{COT1}. They also showed that $\scrf_{2}/\scrf_{2.5}$ has
infinite rank \cite{COT1, COT2}. Later Cochran and Teichner showed
that their filtration is highly nontrivial. That is,
$\scrf_{n}/\scrf_{n.5}$ is infinite for all $n$ \cite{CT}.

In this paper, we give new obstructions for knots being doubly
slice using the ideas of COT. One easily sees that a knot is
doubly slice if and only if there exist two slice disk and 4-ball
pairs whose union along their boundary gives an {\em unknotted}
2-sphere in the 4-sphere. In this regard, for half integers $m$
and $n$, we define a knot to be {\em $(m,n)$-solvable} if the knot
has an $(m)$-solution and an $(n)$-solution such that the union of
these solutions along their boundary gives a closed 4-manifold
whose fundamental group is isomorphic to an infinite cyclic group
(see Definition~\ref{def:bi-solvable}). In particular, we define a
knot to be {\em doubly $(m)$-solvable} if it is $(m,m)$-solvable.
We remark that Freedman \cite{Fr} showed that a 2-knot is
unknotted in the 4-sphere if and only if the fundamental group of
the knot exterior is isomorphic to an infinite cyclic group. We
show that a doubly slice knot is $(m,n)$-solvable for all $m,n$
(Proposition~\ref{prop:obstruction}). For given half-integers
$k\ge m$ and $\ell\ge n$, if a knot is $(k,\ell)$-solvable then it
is $(m,n)$-solvable. This is easily proven since a $(k)$-solution
(respectively an $(\ell)$-solution) for a knot is an
$(m)$-solution (respectively an $(n)$-solution) (refer to
\cite[Remark 1.1.3]{COT1}). Moreover we show that if two knots are
$(m,n)$-solvable, then so is their connected sum
(Proposition~\ref{prop:connected sum}). This implies that if we
denote by $\scrf_{m,n}$ the set of $(m,n)$-solvable knots, then
$\left\{\scrf_{m,n}\right\}_{m,n\ge 0}$ becomes a bi-filtration of
the monoid of knots (under the connected sum operation). We
investigate this bi-filtration and construct examples of knots
showing non-triviality of the bi-filtration at all levels. Our
main theorem is as follows :

\begin{thm}
\label{thm:main}
\begin{enumerate}
\item For a given integer $m\ge 2$, there exists a ribbon knot
(hence slice) $K$ such that $K$ is algebraically doubly slice,
doubly $(m)$-solvable, but not doubly $(m.5)$-solvable. \item For
given integers $k,\ell\ge 2$, there exists an algebraically doubly
slice knot $K$ such that $K$ is $(k,\ell)$-solvable, but neither
$(k.5,\ell)$-solvable, nor $(k,\ell.5)$-solvable.
\end{enumerate}
\end{thm}

\noindent A knot is called a \emph{ribbon knot} if it bounds an
immersed 2-disk (called \emph{ribbon} or \emph{ribbon disk}) in
the 3-sphere with only ribbon singularities. (We say an immersed
2-disk $f(D^2)$ where $f : D^2 \rightarrow S^3$ is an immersion
has \emph{ribbon singularities} if the inverse image of the
singularities consists of pairs of arcs on $D^2$ such that one arc
of each pair is interior to $D^2$.) Note that a ribbon knot is a
slice knot. To see this, push the singular parts of the ribbon
disk into $B^4$ to get a slice disk.

Classical invariants are reflected at the initial levels of the
bi-filtration. In particular, we show that if a knot is doubly
(1)-solvable, then its Blanchfield form is hyperbolic
(Proposition~\ref{prop:1.0}). (It is unknown to the author if the
converse is true.) We also show that a knot has vanishing Arf
invariant if and only if it is doubly (0)-solvable, and
algebraically slice if and only if it is doubly (0.5)-solvable
(Corollary~\ref{cor:0 and 0.5-2}).

To prove the main theorem, we construct a fibred doubly slice knot
of genus 2 which will be called {\em the seed knot}. We choose a
trivial link in the 3-sphere that is disjoint from the seed knot
and choose auxiliary Arf invariant zero knots. Then \emph{genetic
modification} is performed on the seed knot via the chosen trivial
link and auxiliary knots to obtain the desired examples of knots.
This genetic modification is the same as the one used in
\cite{COT2, CT} and will be explained in
Section~\ref{sec:genetic}. In fact, in \cite{CT} Cochran and
Teichner make use of genetic modification to construct the
examples of knots which are $(m)$-solvable but not
$(m.5)$-solvable in COT's filtration of the knot concordance
group. In comparison with their examples, to prove
Theorem~\ref{thm:main}(1), our examples need to be slice, hence
$(k)$-solvable for all $k$. Hence a technical difficulty arises,
and we perform genetic modification in a more sophisticated way
than in \cite{CT}.

To show that a knot is not doubly $(m.5)$-solvable, we use von
Neumann $\rho$-invariants defined by Cheeger and Gromov \cite{CG}.
In particular, we make use of the fact that there is a universal
bound for von Neumann $\rho$-invariants for a fixed 3-manifold
\cite{CT}\cite[Theorem 3.1.1]{R}. More details about this can be
found in Section~\ref{sec:not m.5} and Section~\ref{sec:proof}.

This bi-filtration of knots induces a bi-filtration of the double
concordance group. Two knots $K_1$ and $K_2$ are called {\em
doubly concordant} if $K_1\# J_1$ is isotopic to $K_2\# J_2$ for
some doubly slice knots $J_1$ and $J_2$. (Here `$\#$' means the
connected sum.) This is an equivalence relation, and the
equivalence classes with the connected sum operation form {\em the
double concordance group}. We denote the set of the equivalence
classes represented by $(m,n)$-solvable knots by
$\overline{\scrf}_{m,n}$. We show that each
$\overline{\scrf}_{m,n}$ is a subgroup of the double concordance
group and $\{\overline{\scrf}_{m,n}\}_{m,n\ge 0}$ is a
bi-filtration of the double concordance group
(Corollary~\ref{cor:bi-filtration}).

This paper is organized as follows. In
Section~\ref{sec:definition}, we define $(m,n)$-solvable knots and
show that doubly slice knots are $(m,n)$-solvable for all $m$ and
$n$. We induce a bi-filtration of the monoid of knots and
investigate properties of the bi-filtration at the initial levels.
In Section~\ref{sec:genetic}, we explain how to construct
$(m,n)$-solvable knots using genetic modification. In
Section~\ref{sec:not m.5}, we explain Cochran and Teichner's work
in \cite{CT} and show when $(m)$-solutions are not
$(m.5)$-solutions. In Section~\ref{sec:proof}, we give a proof of
Theorem~\ref{thm:main}. In Section~\ref{sec:bi-filtration}, we
construct a bi-filtration of the double concordance group. Finally
in Section~\ref{sec:1.0} we give the examples of knots
demonstrating non-triviality of the bi-filtration of the monoid of
knots at lower levels.
\\

{\bf Notation.} Throughout this paper, $M_K$ denotes 0-surgery on
a knot $K$ in $S^3$ and $\L$ (respectively $\L'$) denotes the
group ring $\bbz[t,t^{-1}]$ (respectively $\bbq[t,t^{-1}]$). The
set of non-negative integers is denoted by $\bbn_0$. For
convenience we use the same notations for a simple closed curve
and the homotopy (and homology) class represented by the curve.
The integer coefficients are understood for homology groups unless
specified otherwise

\section{$(m,n)$-solvable knots and the basic properties}
\label{sec:definition} $(m,n)$-solvable knots and doubly
$(m)$-solvable knots are defined as follows.
\begin{definition}
\label{def:bi-solvable} Let $m,n\in \frac12\bbn_0$. A 3-manifold
$M$ is called \emph{$(m,n)$-solvable via $(W_1,W_2)$} if $M$ is
$(m)$-solvable via $W_1$ and $(n)$-solvable via $W_2$ such that
the fundamental group of the union of $W_1$ and $W_2$ along their
boundary $M$ is isomorphic to $\bbz$. (i.e., $\pi_1(W_1\cup_{M}
W_2) \cong \bbz$.) A knot $K$ is called \emph{$(m,n)$-solvable via
$(W_1,W_2)$} if $M_K$ is $(m,n)$-solvable via $(W_1,W_2)$. The
ordered pair $(W_1,W_2)$ is called an \emph{$(m,n)$-solution for
$K$} (or $M_K$). The set of all $(m,n)$-solvable knots is denoted
by $\scrf_{m,n}$.
\end{definition}

\begin{definition}
\label{def:doubly solvable} A knot $K$ is \emph{doubly
$(m)$-solvable} if it is $(m,m)$-solvable. An $(m,m)$-solution for
$K$ is called a {\em double $(m)$-solution for $K$}.
\end{definition}

For the reader's convenience, the definition of $(n)$-solvability
is given below. For the related terminologies and more
explanations about $(n)$-solvable knots, refer to \cite{COT1}.

\begin{definition}[COT1]
\label{def:(n)-solvability} Let $n\in \bbn_0$. A $3$-manifold $M$
is \emph{$(n)$-solvable} (resp. \emph{$(n.5)$-solvable}) if there
is an $H_1$-bordism $W$ which contains an $(n)$-Lagrangian (resp.
$(n+1)$-Lagrangian) with $(n)$-duals. If $M$ is zero surgery on a
knot or a link then the corresponding knot or link is called
$(n)$-solvable (resp. $(n.5)$-solvable).
\end{definition}

\begin{remark}
\label{rem:basic}
\begin{itemize}
\item [(i)]By van Kampen Theorem, the condition
$$\pi_1(W_1 \cup_{M} W_2)\cong \bbz$$
is equivalent to the condition that
the following diagram is a push-out diagram in the category of
groups and homomorphisms. In the diagram, $i_1$ and $i_2$ are the
homomorphisms induced from the inclusion maps from $M_K$ into
$W_1$ and $W_2$, and $j_1$ and $j_2$ are the abelianization.
$$
\begin{diagram}\dgARROWLENGTH=1.0em
\node[2]{\pi_1(W_1)} \arrow{se,t}{j_1}\\
\node{\pi_1(M)} \arrow{ne,t}{i_1} \arrow{se,b}{i_2} \node[2]{\bbz}
\\
\node[2]{\pi_1(W_2)} \arrow{ne,b}{j_2}
\end{diagram}
$$
In other words, the condition is equivalent to the condition
$$\pi_1(W_1) \ast_{\pi_1(M)} \pi_1(W_2) \cong \bbz.$$

\item [(ii)] Let $E_K$ be the exterior of $K$ in $S^3$ (i.e.,
$E_K \equiv S^3\setminus N(K)$ where $N(K)$ is an open tubular
neighborhood of $K$). Then $\pi_1(W_1\cup_{M_K} W_2) \cong
\pi_1(W_1\cup_{E_K} W_2)$. This is easily proven using the fact
that $\pi_1(M_K) \cong \pi_1(E_K)/\left<\ell\right>$ where
$\left<\ell\right>$ is the subgroup normally generated by the
longitude $\ell$ of $K$.

\item[(iii)] By definition of $(m)$-solvability,  if $(W_1,W_2)$
is an $(m,n)$-solution, $W_1$ and $W_2$ are spin 4-manifolds, and
one easily sees that $W_1 \cup_{M} W_2$ is spinable. But we do not
need this fact for our purpose.

\item[(iv)] If a knot $K$ is $(m,n)$-solvable, then one easily
sees that $K$ is $(k)$-solvable where $k$ is the maximum of $m$
and $n$.
\end{itemize}
\end{remark}
The following proposition shows that doubly slice knots are
contained in the intersection of all $\scrf_{m,n}$'s.
\begin{prop}
\label{prop:obstruction}If a knot $K$ is doubly slice, then it is
$(m,n)$-solvable for all $m$ and $n$.
\end{prop}
\begin{proof}Since $K$ is doubly slice, there are two slice disk
and 4-ball pairs $(B^4_1, D^2_1)$ and $(B^4_2, D^2_2)$ such that
$(S^3,K) = \partial(B^4_1, D^2_1) = \partial(B^4_2, D^2_2)$ and
$D^2_1\cup_K D^2_2$ is an unknotted 2-sphere in the 4-sphere.
Since the second homology of a slice disk exterior is trivial,
every slice disk exterior is an $(m)$-solution for the knot for
all $m$ (see \cite[Remark 1.3.1]{COT1}). So if we let $W_i \equiv
B^4_i\setminus N(D^2_i)$ for $i = 1,2$, then we may think that
$W_1$ is an $(m)$-solution and $W_2$ is an $(n)$-solution for a
given pair of half-integers $m$ and $n$. Furthermore,
$W_1\cup_{E_K}W_2$ is homeomorphic to the exterior of an unknotted
2-sphere in the 4-sphere (which is homeomorphic to $S^1\times
D^3$), hence $(W_1,W_2)$ satisfies the required fundamental group
condition.
\end{proof}
\noindent The following proposition shows that $\scrf_{m,n}$ is a
submonoid of the monoid of knots under the connected sum
operation.

\begin{prop}
\label{prop:connected sum} Suppose $K$ and $J$ are
$(m,n)$-solvable knots. Then $K\#J$ is $(m,n)$-solvable.
\end{prop}
\begin{proof}
Let $(V_1,V_2)$ be an $(m,n)$-solution for $K$ and $(W_1,W_2)$ be
an $(m,n)$-solution for $J$. We will construct a specific
$(m,n)$-solution for $K\#J$ using these solutions. We begin by
constructing a standard cobordism $C$ between $M_K\amalg M_J$ and
$M_{K\#J}$. Start with $(M_K\amalg M_J)\times [0,1]$ and add a
1-handle to $(M_K\amalg M_J) \times \{1\}$ such that the upper
boundary is a connected 3-manifold given by surgery on a split
link $K\amalg J$ with 0-framing. Next, add a 2-handle with
0-framing to the upper boundary along an unknotted circle which
wraps around $K$ and $J$ once. (This equates the meridional
generators of the first homology of $M_K$ and $M_J$.) The
resulting 4-manifold is $C$. That is, $\partial_{-}C = M_K\amalg
M_J$ and $\partial_{+}C = M_{K\#J}$. See \cite[Theorem 4.1]{COT2}
and its proof for more details.

Now let $X_i$ be the union of $C$, $V_i$, and $W_i$ along the
boundaries as shown in Figure~\ref{fig:C} for $i=1,2$. We claim
that $(X_1,X_2)$ is an $(m,n)$-solution for $K\#J$. First we show
that $X_1$ is an $(m)$-solution for $K\#J$. (The proof that $X_2$
is an $(n)$-solution for $K\#J$ will follow similarly.) In the
construction of the cobordism $C$, one can see that $H_1(C)\cong
\bbz$ and the inclusion from any boundary component of $C$ induces
an isomorphism. It follows that the inclusion induced map
$H_1(M_{K\#J}) \rightarrow H_1(X_1)$ is an isomorphism. Since
adding a 1-handle and a 2-handle has no effect on $H_2$,
$H_2(C)\cong H_2(M_K) \oplus H_2(M_J)$. Let $Y_1 \equiv V_1 \amalg
W_1$. From the pair of spaces $(C, Y_1)$, we get the following
Mayer-Vietoris sequence
$$
\cdots \rightarrow H_2(M_K \amalg M_J) \rightarrow H_2(C) \oplus
H_2(Y_1) \rightarrow H_2(X_1) \rightarrow H_1(M_K \amalg M_J)
\rightarrow \cdots .
$$
Since $H_1(M_K \amalg M_J) \rightarrow H_1(Y_1)$ is an
isomorphism, $H_2(X_1) \rightarrow H_1(M_K \amalg M_J)$ is the
zero map. By the above observation on $H_2(C)$, $H_2(M_K \amalg
M_J) \rightarrow H_2(C)$ is an isomorphism, hence surjective.
Since the boundary map $H_3(V_1, M_K) \rightarrow H_2(M_K)$ is the
dual of an isomorphism $H^1(V_1) \rightarrow H^1(M_K)$, it is an
isomorphism. Hence $H_2(M_K) \rightarrow H_2(V_1)$ is the zero
map. Similarly $H_2(M_J) \rightarrow H_2(W_1)$ is the zero map,
thus so is $H_2(M_K \amalg M_J) \rightarrow H_2(Y_1)$. So
$H_2(X_1) \cong H_2(Y_1) \cong H_2(V_1) \oplus H_2(W_1)$. Since
the intersection form on $Y_1$ splits naturally on $V_1$ and
$W_1$, the ``union" of the $(m)$-Lagrangians and $(m)$-duals for
$V_1$ and $W_1$ forms the $(m)$-Lagrangian and $(m)$-dual for
$X_1$. So $X_1$ is an $(m)$-solution for $K\#J$. (For more details
on $(m)$-solutions, $(m)$-Lagrangians, and $(m)$-duals, refer to
\cite[Section 7,8]{COT1}.)

\begin{figure}[htb]
\includegraphics[scale=1]{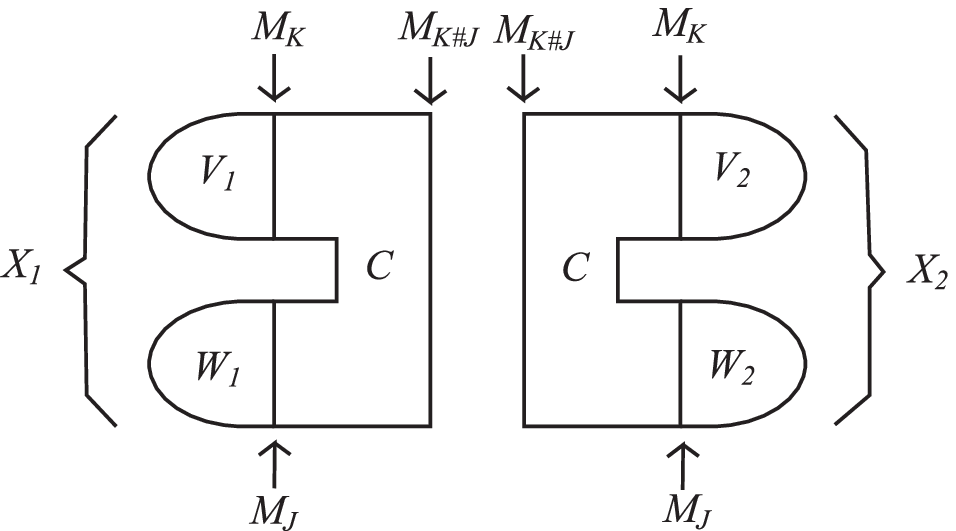}
\caption{}\label{fig:C}
\end{figure}

It remains to show $\pi_1(X_1\cup_{M_{K\#J}} X_2) \cong \bbz$, or
$\pi_1(X_1\cup_{E_{K\#J}} X_2) \cong \bbz$ by
Remark~\ref{rem:basic}. Thus to prove the proposition, it is
enough to show that the following diagram is a push-out diagram.
$$
\begin{diagram}\dgARROWLENGTH=1.0em
\node{\pi_1(E_{K\#J})} \arrow{e,t}{i_1} \arrow{s,l}{i_2}
\node{\pi_1(X_1)} \arrow{s,r}{j_1}
\\
\node{\pi_1(X_2)} \arrow{e,b}{j_2} \node{\bbz}
\end{diagram}
$$

Since $(V_1,V_2)$ and $(W_1,W_2)$ are $(m,n)$-solutions for $K$
and $J$ respectively, we have the following push-out diagrams.
$$
\begin{diagram}\dgARROWLENGTH=1.0em
\node{\pi_1(E_{K})} \arrow{e} \arrow{s} \node{\pi_1(V_1)}
\arrow{s} \node{\pi_1(E_{J})} \arrow{e} \arrow{s}
\node{\pi_1(W_1)} \arrow{s}
\\
\node{\pi_1(V_2)} \arrow{e} \node{\bbz} \node{\pi_1(W_2)}
\arrow{e} \node{\bbz}
\end{diagram}
$$

\noindent By taking free products and factoring out each group by
the normal subgroup $\left<\mu_K\mu_J^{-1}\right>$, we have the
following push-out diagram. (Here $\mu_K$ and $\mu_J$ are
meridians of $K$ and $J$ respectively.)
$$
\begin{diagram}\dgARROWLENGTH=1.0em
\node{\frac{\pi_1(E_K)\ast\pi_1(E_J)}{\left<\mu_K
\mu_J^{-1}\right>}} \arrow{e} \arrow{s}
\node{\frac{\pi_1(V_1)\ast\pi_1(W_1)}{\left<\mu_K
\mu_J^{-1}\right>}} \arrow{s}
\\
\node{\frac{\pi_1(V_2)\ast\pi_1(W_2)}{\left<\mu_K
\mu_J^{-1}\right>}} \arrow{e} \node{\bbz}
\end{diagram}
$$
\noindent Note that
$$
\frac{\pi_1(E_K)\ast\pi_1(E_J)}{\left<\mu_K \mu_J^{-1}\right>}
\cong \pi_1(E_{K\#J}).
$$
By the construction of the cobordism $C$,
$$
\frac{\pi_1(V_i)\ast\pi_1(W_i)}{\left<\mu_K \mu_J^{-1}\right>}
\cong \pi_1(X_i)
$$
for $i=1,2$.

\end{proof}
\noindent From Proposition~\ref{prop:obstruction} and
Proposition~\ref{prop:connected sum}, we can easily deduce the
following corollary.
\begin{cor}
\label{cor:bi-filtration of knots}The family
$\{\scrf_{m,n}\}_{m,n\ge 0}$ is a bi-filtration of the monoid of
knots under the connected sum operation where doubly slice knots
lie in the intersection of all $\scrf_{m,n}$'s.
\end{cor}

Next, we study the properties of this bi-filtration at lower
levels.
\begin{prop}
\label{prop:0 and 0.5} Suppose $n = 0$ or $0.5$. Then a knot $K$
is doubly $(n)$-solvable if and only if it is $(n)$-solvable.
\end{prop}
\begin{proof}
One direction is clear by Remark~\ref{rem:basic}\,(iv). For the
other direction, suppose $K$ is $(n)$-solvable via $W$. By doing
surgery on the commutator subgroup of $\pi_1(W)$ (note that the
commutator subgroup is finitely normally generated), we may assume
that $\pi_1(W) \cong \bbz$. Let $W_1$ and $W_2$ be copies of $W$.
Their fundamental groups are isomorphic to $\bbz$ and generated by
the meridian of $K$. So using van Kampen Theorem, one sees that
$\pi_1(W_1\cup_{M_K} W_2) \cong\bbz$, hence $K$ is doubly
$(n)$-solvable via $(W_1,W_2)$.
\end{proof}

\noindent It is known that a knot is $(0)$-solvable if and only if
it has vanishing Arf invariant, and $(0.5)$-solvable if and only
if it is algebraically slice (that is, its associated Seifert
forms are metabolic). (See \cite{COT1}.) So we have the following
corollary.

\begin{cor}
\label{cor:0 and 0.5-2} A knot is doubly $(0)$-solvable if and
only if it has vanishing Arf invariant, and doubly
$(0.5)$-solvable if and only if it is algebraically slice.
\end{cor}

We investigate relationship between the bi-filtration
$\{\scrf_{m,n}\}_{m,n\ge 0}$ and algebraically doubly slice knots.
For a knot $K$, we have the (nonsingular and sesquilinear)
Blanchfield form $B\ell : H_1(M_K;\L)\times H_1(M_K;\L)\rightarrow
\bbq(t)/\L$ (see \cite{B}). For a Seifert matrix of $K$, say $A$,
the Blanchfield form is presented by $(1-t)(tA -A^{T})^{-1}$ (see
\cite{L}). Kearton showed that this presentation matrix is doubly
null cobordant if and only if the Seifert matrix $A$ is
$S$-equivalent to a doubly null cobordant Seifert matrix (see
\cite{Ke1} and \cite{Tr}). A matrix is called {\em doubly null
cobordant} if it is congruent by an integer unimodular matrix to a
matrix of the form
$$
\left(\begin{matrix} 0 & *\cr * & 0
\end{matrix}\right).
$$
This implies that the Blanchfield form of a knot $K$ is hyperbolic
(that is, $H_1(M_K;\L) = A \oplus B$ where $A$ and $B$ are
$\L$-submodules of $H_1(M;\L)$ and they are self-annihilating with
respect to $B\ell$) if and only if $K$ has a Seifert matrix which
is S-equivalent to a doubly null cobordant matrix. Now we have the
following proposition.

\begin{prop}
\label{prop:1.0} Suppose a knot $K$ is doubly $(1)$-solvable via
$(W_1,W_2)$. Let $i_j : M_K \rightarrow W_j$ be the inclusion map
for $j= 1,2$. Then
$$
H_1(M_K;\L) = \Ker(i_1)_* \oplus \Ker(i_2)_*.
$$
Furthermore, $\Ker(i_1)_* \cong H_1(W_2;\L)$ and $\Ker(i_2)_*
\cong H_1(W_1;\L)$. Moreover, for each $j$, $\Ker(i_j)_*$ is a
self-annihilating submodule (that is, $\Ker(i_j)_* =
\Ker(i_j)_*^\bot$) with respect to the Blanchfield form
$$
B\ell : H_1(M_K;\L) \times H_1(M_K;\L)\rightarrow \bbq(t)/\L.
$$
Hence the Blanchfield form $B\ell$ is hyperbolic.
\end{prop}
\begin{proof}
Let $W$ be $W_1\cup_{M_K} W_2$. Recall that $\L' \equiv \bbq
[t,t^{-1}]$. Since $W_1\cap W_2 = M_K$, we have the following
Mayer-Vietoris sequence.
$$\cdots \rightarrow H_2(W;\L')\xrightarrow{\partial}
H_1(M_K;\L')\xrightarrow{f} H_1(W_1;\L')\oplus
H_1(W_2;\L')\xrightarrow{g} H_1(W;\L')\rightarrow\cdots$$ Since
$\pi_1(W)\cong\bbz$, $H_1(W;\L') = \{0\}$. We show that $f$ is
injective. Suppose $x\in \Ker f$. We can consider $f$ as
$((i_1)_*,(i_2)_*)$. Therefore $x\in \Ker (i_j)_*$ for $j = 1,2$.
By \cite[Theorem 3.5 and Theorem 3.6]{COT1}, $x$ induces a
representation $\phi : \pi_1(M) \rightarrow \Gamma^U_1$ where
$\Gamma^U_1\equiv (\bbq(t)/\L')\rtimes\bbz$ such that $\phi$ can
be extended to $\Phi_1 : \pi_1(W_1)\rightarrow \Gamma^U_1$ and
$\Phi_2 : \pi_1(W_2)\rightarrow \Gamma^U_1$, hence we have the
following commutative (push-out) diagram.
$$
\begin{diagram}\dgARROWLENGTH=1.5em
\node[2]{\pi_1(W_1)} \arrow{se,b}{j_1} \arrow{ese,t}{\Phi_1}\\
\node{\pi_1(M_K)} \arrow{ne,t}{i_1} \arrow{se,b}{i_2}
\node[2]{\bbz} \arrow{e,t}{\alpha} \node{\Gamma^U_1}\\
\node[2]{\pi_1(W_2)} \arrow{ne,t}{j_2} \arrow{ene,b}{\Phi_2}
\end{diagram}
$$
In this diagram, we get the homomorphism $\alpha$ by the universal
property of the push-out diagram. Let $\e : \pi_1(M_K)
\rightarrow\bbz$ be the abelianization (in fact, $\e = j_1\circ
i_1$). For $y\in\pi_1(M_K)$, $\phi(y)$ is calculated as $\phi(y) =
(B\ell'(x,y\mu^{-\e(y)}), \e(y))$ for a meridian $\mu$ of the knot
$K$ and the {\em rational} Blanchfield pairing
$$
B\ell' : H_1(M_K;\L')\times H_1(M_K;\L')\rightarrow \bbq(t)/\L'.
$$
Thus $\phi(\mu) = (0,1)\in\Gamma^U_1$. By the commutativity of the
diagram, we have $\alpha(1) = (0,1) \in \Gamma^U_1$. Thus for any
meridian, say $\mu'$, of the knot, $\phi(\mu') = \alpha(1) = (0,
1)$ in $\Gamma^U_1$. Thus $\Im$ $\alpha = \{0\}\rtimes \bbz
\subset \Gamma^U_1$, hence $\phi(y)\in \{0\}\rtimes \bbz$ for all
$y\in \pi_1(M_K)$. Therefore $B\ell'(x,x') = 0$ for all $x'\in
H_1(M_K;\L')$. Since the rational Blanchfield pairing is
nonsingular, this implies $x=0$, hence $f$ is injective. Hence
$H_1(M_K;\L') = \Ker (i_1)_* \oplus \Ker (i_2)_*$ where $\Ker
(i_1)_*\cong H_1(W_2;\L')$ and $\Ker (i_2)_* \cong H_1(W_1;\L')$.

Now we replace the coefficients $\L'$ by $\L$. One sees that
$H_1(W;\L) = \{0\}$ because $\pi_1(W)\cong\bbz$. The homomorphism
$f$ is still injective since $H_1(M;\L)$ is $\bbz$-torsion free.
Therefore $H_1(M_K;\L) = \Ker (i_1)_* \oplus \Ker (i_2)_*$ where
$\Ker (i_1)_* \cong H_1(W_2;\L)$, $\Ker (i_2)_* \cong
H_1(W_1;\L)$. We need to show that $\Ker (i_j)_*$ is
self-annihilating for each $j$. Since $W_j$ is an (integral)
$(1)$-solution for $K$,
$$
TH_2(W_j,M_K;\L) \xrightarrow{\partial} H_1(M_K;\L)
\xrightarrow{(i_j)_*} H_1(W_j;\L)
$$
is exact by \cite[Lemma 4.5]{COT1} where $TH_2$ denotes the
$\L$-torsion submodule. Note that the Kronecker map
$$
\kappa :
H^1(W_j;\bbq(t)/\L) \rightarrow \Hom_\L(H_1(W_j;\L),\bbq(t)/\L)
$$
is an isomorphism from the universal coefficient spectral sequence
and the map
$$
(i_j)^\# : \Hom_\L(H_1(W_j;\L), \bbq(t)/\L)
\rightarrow \Hom_\L(H_1(M;\L)/\Ker(i_j)_*, \bbq(t)/\L)
$$
is also an isomorphism since $(i_j)_* : H_1(M;\L) \rightarrow
H_1(W_j;\L)$ is onto. Now one follows the course of the proof of
\cite[Theorem 4.5]{COT1} and obtains that $\Ker(i_j)_* =
(\Ker(i_j)_*)^\bot$.

\end{proof}

By the observation preceding Proposition~\ref{prop:1.0}, we have
the following corollary.

\begin{cor}
\label{cor:1.0} If a knot $K$ is doubly $(1)$-solvable, then $K$
has a Seifert matrix which is S-equivalent to a doubly null
cobordant matrix.
\end{cor}
It is unknown to the author if a knot with the hyperbolic
Blanchfield form is doubly $(1)$-solvable.

\begin{remark}
\label{rem:s-equiv} That a matrix is $S$-equivalent to a doubly
null cobordant matrix does not imply that the matrix itself is
doubly null cobordant. Thus that a knot is algebraically doubly
slice does not mean that all of its associated Seifert forms are
hyperbolic (but at least there is one Seifert form that is
hyperbolic). (See \cite{Ke1}.)
\end{remark}

\section{genetic modification}
\label{sec:genetic} In this section we recall the notion of
\emph{genetic modification} and show when it preserves
$(m,n)$-solvability of a knot. This modification of knots is the
same as the one used in \cite{COT2} and \cite{CT}.

Let $K$ be a knot in $S^3$. Let $\eta$ be a trivial knot in $S^3$
which is disjoint from $K$. Let $J$ be another knot. Take the
exterior of $\eta$ (which is homeomorphic to a solid torus) and
the exterior of $J$. Now identify them along their boundary such
that the meridian of $\eta$ (say $\mu_\eta$) is identified with
the longitude of $J$ (say $\ell_J$) and the longitude of $\eta$
(say $\ell_\eta$) is identified with the meridian of $J$ (say
$\mu_J$). The resulting ambient manifold is homeomorphic to $S^3$,
and we denote the image of $K$ under this modification by $K(J,
\eta)$. In fact, $K(J,\eta)$ is a satellite of $J$. This
construction can be generalized to the case that we have a trivial
link $\{\eta_1, \eta_2, \ldots, \eta_n\}$ which misses $K$ and a
set of auxiliary knots $\{J_1, J_2, \ldots, J_n\}$ by repeating
the construction. We denote the resulting knot by $K(\{J_1, J_2,
\ldots, J_n\}, \{\eta_1, \eta_2, \ldots, \eta_n\})$. More details
can be found in \cite{COT2}.

The following proposition is implicitly proved in \cite{COT2}. For
a group $G$, we define $G^{(0)}\equiv [G,G]$, and inductively
$G^{(n+1)}\equiv [G^{(n)}, G^{(n)}]$ for $n\ge 0$. That is,
$G^{(n)}$ is {\em the $n$-th derived subgroup} of $G$.
\begin{prop}\cite[Propositin 3.1]{COT2}
\label{prop:(n)-solvable} If $K$ is $(n)$-solvable via $W$,
$\eta\in\pi_1(W)^{(n)}$, and $J$ is a knot with vanishing Arf
invariant, then $K(J, \eta)$ is $(n)$-solvable.
\end{prop}

\noindent We give a brief explanation as to how to construct an
$(n)$-solution for $K(J, \eta)$ from $W$ in the above proposition.
This will also serve to set the notations that will be used later
in this paper. Since Arf invariant vanishes for $J$, $J$ is
$(0)$-solvable. Let $W_J$ be a $(0)$-solution for $J$. By doing
surgery on the commutator subgroup of $\pi_1(W_J)$, we may assume
that $\pi_1(W_J)\cong\bbz$. Note that $\partial W = M_K$ and
$\partial W_J = M_J = E_J\cup S^1\times D^2$ where $E_J$ is the
exterior of $J$, $\{\ast\}\times \partial D^2$ is the longitude
$\ell_J$, and $S^1\times \{\ast\}$ is the meridian $\mu_J$. Let
$\eta\times D^2$ be a tubular neighborhood of $\eta$ in $M_K$.
Then the $(n)$-solution for $K(J, \eta)$, say $W'$, is obtained
from $W$ and $W_J$ by identifying $\eta\times D^2\subset
\partial W$ and $S^1\times D^2\subset
\partial W_J$.

The next proposition shows that we have a similar result for
$(m,n)$-solvable knots. In the statement, $W_1'$ and $W_2'$ denote
the $(m)$-solution and the $(n)$-solution for $K(J,\eta)$ obtained
from $W_1$ and $W_2$ by the above construction in the previous
paragraph.
\begin{prop}
\label{prop:bi-solvable} Suppose $K$ is $(m,n)$-solvable via
$(W_1,W_2)$, $\eta\in\pi_1(W_1)^{(m)}\cap \pi_1(W_2)^{(n)}$, and
$J$ is a knot with vanishing Arf invariant. Then $K'=K(J,\eta)$ is
$(m,n)$-solvable via $(W_1',W_2')$.
\end{prop}
\begin{proof}By Proposition~\ref{prop:(n)-solvable}, $W_1'$ and
$W_2'$ are an $(m)$-solution and an $(n)$-solution for $K'$,
respectively. Let $W \equiv W_1'\cup_{M_{K'}}W_2'$. We need to
show that $\pi_1(W)\cong\bbz$. For convenience, let $M \equiv M_K$
and $M' \equiv M_{K'}$. Since $K$ is $(m,n)$-solvable via
$(W_1,W_2)$, we have the following push-out diagram in the
category of groups and homomorphisms.
$$
\begin{diagram}\dgARROWLENGTH=1.0em
\node[2]{\pi_1(W_1)} \arrow{se,t}{j_1}\\
\node{\pi_1(M)} \arrow{ne,t}{i_1} \arrow{se,b}{i_2} \node[2]{\bbz}
\\
\node[2]{\pi_1(W_2)} \arrow{ne,b}{j_2}
\end{diagram}
$$
We will show that the following diagram is also a push-out
diagram, then this will complete the proof. In the diagram, $i_1'$
and $i_2'$ are the homomorphisms induced from the inclusions and
$j'_1$ and $j'_2$ are the abelianization.
$$
\begin{diagram}\dgARROWLENGTH=1.0em
\node[2]{\pi_1(W_1')} \arrow{se,t}{j_1'}\\
\node{\pi_1(M')} \arrow{ne,t}{i_1'} \arrow{se,b}{i_2'}
\node[2]{\bbz}
\\
\node[2]{\pi_1(W_2')} \arrow{ne,b}{j_2'}
\end{diagram}
$$
Suppose we are given a commutative diagram as below where $\Gamma$
is a group.
$$
\begin{diagram}\dgARROWLENGTH=1.0em
\node[2]{\pi_1(W_1')} \arrow{se,t}{\alpha_1}\\
\node{\pi_1(M')} \arrow{ne,t}{i_1'} \arrow{se,b}{i_2'}
\node[2]{\Gamma}
\\
\node[2]{\pi_1(W_2')} \arrow{ne,b}{\alpha_2}
\end{diagram}
$$
We study relationship among the fundamental groups of the spaces.
Observe that $M' = (M \setminus \int(\eta\times
D^2))\cup_{\eta\times S^1} E_J$ where $\eta\times S^1 =
\partial (\eta\times D^2)$. Let $X = M\setminus\int(\eta\times
D^2)$. By van Kampen Theorem, $\pi_1(M) \cong
\pi_1(X)/\left<\mu_\eta\right>$ where $\left<\mu_\eta\right>$ is
the subgroup normally generated by $\mu_\eta$ in $\pi_1(X)$, and
$$
\pi_1(M')\cong
\displaystyle\frac{\pi_1(X)\ast\pi_1(E_J)}{\left<\ell_\eta\mu_J^{-1},
\mu_\eta\ell_J^{-1}\right>}.
$$
For $W_1'$ and $W_2'$, van Kampen Theorem shows that for $i=1,2$,
$$
\pi_1(W_i')\cong \displaystyle\frac{\pi_1(W_i)\ast\pi_1(W_J)}
{\left<\ell_\eta\mu_J^{-1}\right>} \cong \displaystyle
\frac{\pi_1(W_i)\ast\left<\mu_J\right>}
{\left<\ell_\eta\mu_J^{-1}\right>} \cong \pi_1(W_i).
$$
For simplicity, let
$$
G \equiv \displaystyle
\frac{\pi_1(X)\ast\pi_1(E_J)}{\left<\ell_\eta\mu_J^{-1},
\mu_\eta\ell_J^{-1}\right>}
$$
and $f : G\rightarrow \pi_1(M')$ be the isomorphism given by van
Kampen Theorem. Consider the following commutative diagram.
$$
\begin{diagram}\dgARROWLENGTH=1.0em
\node[2]{\pi_1(W_1')} \arrow{se,t}{\alpha_1}\\
\node{G} \arrow{ne,t}{i_1'\circ f} \arrow{se,b}{i_2'\circ f}
\node[2]{\Gamma}
\\
\node[2]{\pi_1(W_2')} \arrow{ne,b}{\alpha_2}
\end{diagram}
$$
Since $\ell_J = e$ in $\pi_1(W_J)$, $\ell_J = e$ in $\pi_1(W_1')$
and $\pi_1(W_2')$. Furthermore, $\pi_1(E_J)$ is mapped into
$\left<\mu_J\right>$ ($= \pi_1(W_J)$) in $\pi_1(W_i')$. Thus
$i_1'\circ f$ and $i_2'\circ f$ factor through
$\displaystyle\frac{\pi_1(M)\ast\left<\mu_J\right>}
{\left<\ell_\eta\mu_J^{-1}\right>}$ which is isomorphic to
$\pi_1(M)$. So we have the following commutative diagram.
$$
\begin{diagram}\dgARROWLENGTH=1.0em
\node[3]{\pi_1(W_1')} \arrow{se,t}{\alpha_1}\\
\node{G} \arrow{ene,t}{i_1'\circ f} \arrow{ese,b}{i_2'\circ f}
\arrow{e} \node{\pi_1(M)} \arrow{ne,b}{k_1} \arrow{se,t}{k_2}
\node[2]{\Gamma}\\
\node[3]{\pi_1(W_2')} \arrow{ne,b}{\alpha_2}
\end{diagram}
$$
Let $p_1 : \pi_1(W_1')\rightarrow\pi_1(W_1)$ be the inverse of the
isomorphism $\pi_1(W_1) \rightarrow \pi_1(W_1')$ induced from the
inclusion. Define $p_2 : \pi_1(W_2')\rightarrow\pi_1(W_2)$
similarly. Then the above diagram induces the following
commutative diagram.
$$
\begin{diagram}\dgARROWLENGTH=1.0em
\node[2]{\pi_1(W_1)} \arrow{se,t}{\alpha_1\circ p_1^{-1}}\\
\node{\pi_1(M)} \arrow{ne,t}{p_1\circ k_1} \arrow{se,b}{p_2\circ
k_2}
\node[2]{\Gamma}\\
\node[2]{\pi_1(W_2)} \arrow{ne,b}{\alpha_2\circ p_2^{-1}}
\end{diagram}
$$
One sees that $i_1 =p_1\circ k_1$ and $i_2 =p_2\circ k_2$. By the
universal property of the push-out diagram, we have a unique
homomorphism $\beta : \bbz\rightarrow\Gamma$ that makes the
following diagram commutative.
$$
\begin{diagram}\dgARROWLENGTH=1.0em
\dgLABELOFFSET=.20ex \dgARROWPARTS=4
\node[2]{\pi_1(W_1)} \arrow{se,b}{j_1} \arrow{ese,t}{\alpha_1\circ p_1^{-1}}\\
\node{\pi_1(M)} \arrow{ne,t}{i_1} \arrow{se,b}{i_2}
\node[2]{\bbz} \arrow{e,t,1}{\beta} \node{\Gamma}\\
\node[2]{\pi_1(W_2)} \arrow{ne,t}{j_2} \arrow{ene,b}{\alpha_2\circ
p_2^{-1}}
\end{diagram}
$$
Thus the following diagram is also commutative
$$
\begin{diagram}\dgARROWLENGTH=1.0em
\dgLABELOFFSET=.20ex \dgARROWPARTS=4
\node[2]{\pi_1(W_1')} \arrow{se,b}{j_1'} \arrow{ese,t}{\alpha_1}\\
\node{\pi_1(M')} \arrow{ne,t}{i_1'} \arrow{se,b}{i_2'}
\node[2]{\bbz} \arrow{e,t,1}{\beta} \node{\Gamma}\\
\node[2]{\pi_1(W_2')} \arrow{ne,t}{j_2'} \arrow{ene,b}{\alpha_2}
\end{diagram}
$$
where $j_1'\equiv j_1\circ p_1$, $j_2'\equiv j_2\circ p_2$. The
choice of $\beta : \bbz\rightarrow\Gamma$ is unique because it is
unique in the previous diagram involving $\pi_1(M)$, $\pi_1(W_1)$,
and $\pi_1(W_2)$.
\end{proof}

\noindent Note that $\pi_1(W_i') \cong \pi_1(W_i)$ in the above
proof. Therefore by applying Proposition~\ref{prop:bi-solvable}
repeatedly, we obtain the following corollary.
\begin{cor}
\label{cor:bi-solvable} Suppose $K$ is $(m,n)$-solvable via
$(W_1,W_2)$. Suppose $\eta_i\in\pi_1(W_1)^{(m)} \cap
\pi_1(W_2)^{(n)}$, and the Arf invariant vanishes for $J_i$ for
$1\le i\le n$. Then
$$
K(\{J_1,J_2,\ldots,
J_n\},\{\eta_1,\eta_2,\ldots, \eta_n\})
$$
is $(m,n)$-solvable via
$(W_1',W_2')$.
\end{cor}

The following lemma and proposition give conditions under which
the knot resulting from genetic modification performed on a ribbon
knot is still a ribbon knot. Let $f_i : D^2 \rightarrow S^3$ be
immersions, $1 \le i \le n$, where each immersed disk $f_i(D^2)$
has only ribbon singularities. We say $f_i(D^2)$ have \emph{ribbon
intersections} if $f_i^{-1}\left(f_i(D^2) \cap f_j(D^2)\right)$,
$i\ne j$, consists of arcs on $D^2$ either having endpoints on
$\partial D^2$ or interior to $D^2$. Recall that $\eta_i$, $1\le i
\le n$, denotes a trivial link which misses a knot $K$. Let $J_i$,
$1 \le i \le n$, denote knots in $S^3$ (not necessarily with
vanishing Arf invariant), and $B^4$ denote the standard 4-ball.

\begin{lem}
\label{lem:ribbon} Suppose $K$ is a ribbon knot bounding a ribbon
disk $B$. Let $\eta_i$, $1 \le i \le n$, bound disjoint embedded
disks $D_i$ in $S^3$ such that $D_i$ and $B$ have ribbon
intersections. Let $J_i$, $1 \le i \le n$, be knots in $S^3$. Then
$K' \equiv K\left(\{J_i\}_{1\le i \le n}, \{\eta_i\}_{1\le i \le
n}\right)$ is a ribbon knot.
\end{lem}
\begin{proof}
Since $D_i$ and $B$ have ribbon intersections, a component of the
intersection of $D_i$ with $B$ is an arc on $D_i$ either having
end points on $\partial D_i$, say \emph{a type I arc}, or interior
to $D_i$, say \emph{a type II arc}. We claim that we may assume
the intersection of $D_i$ with $B$ is only type II arcs. We use an
``outermost arc argument" to show this. Denote type I intersection
arcs of $D_i$ with $B$ by $\alpha_j$, $1 \le j \le \ell$. Suppose
$\alpha_1$ is an outermost arc. That is, $\alpha_1$ splits $D_i$
into two disks, say $A_1$ and $A_2$, such that $A_1$ intersects
the ribbon disk $B$ in only type II arcs. See
Figure~\ref{fig:ribbon} below.

\begin{figure}[htb]
\includegraphics[scale=0.9]{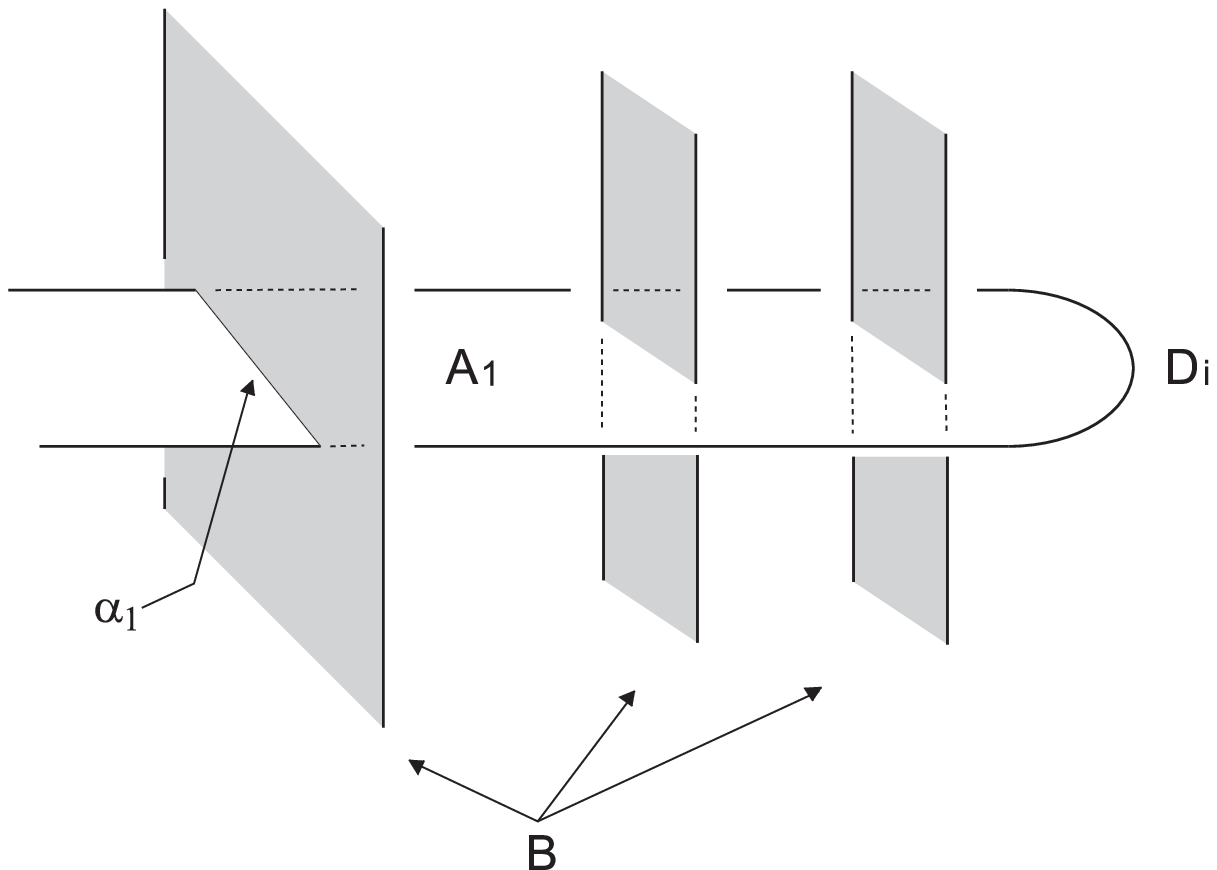}
\caption{}\label{fig:ribbon}
\end{figure}

\noindent Now deform the interior of $B$ along $A_1$ using a
finger move and remove the intersection arc $\alpha_1$. This may
introduce new self-intersections for $B$. But since the
intersection of $A_1$ with $B$ consists of only type II arcs, the
new self-intersections for $B$ are ribbon singularities. Hence the
deformed (immersed) disk is a ribbon disk and it has the same
boundary $K$ as $B$. We repeat this process until we remove all
type I intersection arcs on $D_i$ and this proves the claim.

Observe that $K'$ is indeed the result of cutting open $K$ along
$D_i$ and tying all the strands that pass through $D_i$ into $J_i$
with 0-frame. By the claim the intersection of $D_i$ with $B$ is
arcs interior to $D_i$, hence $B$ passes through the interior of
$D_i$ like bands. (See the two bands on the right in
Figure~\ref{fig:ribbon}.) Thus by cutting open $B$ along $D_i$ and
tying the bands into $J_i$, we obtain an immersed disk, say $B'$,
which is bounded by $K'$. One easily sees that tying $B$ into
$J_i$ does not introduce new self-intersections. Hence $B'$ is
still a ribbon disk.
\end{proof}

The following proposition is due to Peter Teichner.
\begin{prop}[Teichner]
\label{prop:ribbon} Suppose $K$ is a ribbon knot bounding a ribbon
disk $B$. Let $\tilde{B}$ be a slice disk for $K$ obtained by
deforming the ribbon disk $B$ into $B^4$. Suppose $\eta_i$ ($1\le
i \le n$) are knots in $S^3\setminus K$ that are homotopically
trivial in $B^4 \setminus \tilde{B}$. Then there exists a trivial
link $\tau_i$ $(1\le i \le n)$ in $S^3$ which is disjoint from $K$
such that each $\tau_i$ is homotopic to $\eta_i$ in $S^3 \setminus
K$ and $K\left(\{J_i\}_{1\le i \le n}, \{\tau_i\}_{1\le i \le
n}\right)$ is a ribbon knot.
\end{prop}
\begin{proof}
We may think of $B$ as an (immersed) band sum of embedded disks in
$S^3$. Note that the inclusion induced homomorphism
$\pi_1(S^3\setminus K) \rightarrow \pi_1(B^4 \setminus \tilde{B})$
has the kernel that is normally generated by the meridians to the
bands of $B$. Hence there is a trivial link $\tau_i$ in $S^3$
which is disjoint from $K$ such that each $\tau_i$ is homotopic to
$\eta_i$ in $S^3 \setminus K$ and $\tau_i$ ($1\le i \le n$) bound
mutually disjoint embedded disks, say $D_i$, in $S^3$ where each
$D_i$ is obtained by taking a band sum of copies of the meridional
disks to the bands of $B$. One sees that $D_i$ and $B$ have ribbon
intersection. Now the proposition follows from
Lemma~\ref{lem:ribbon}.
\end{proof}

\section{$(m)$-solutions that are not $(m.5)$-solutions}
\label{sec:not m.5} Throughout this section, we assume $K$ is a
genus 2 fibred knot that is $(m)$-solvable. Let
$$
K'\equiv K(\{J_i\}_{1\le i\le n}, \{\eta_i\}_{1\le i\le n}),
$$
the knot resulting from genetic modification. We assume that all
$J_i$ are $(0)$-solvable and $\eta_i$ are lying in
$\pi_1(M_K)^{(m)}$. By Proposition~\ref{prop:(n)-solvable}, $K'$
is $(m)$-solvable. Let $V$ be an $(m)$-solution for $K'$. In this
section we investigate conditions under which it is guaranteed
that $V$ is not an $(m.5)$-solution for $K'$. The key result is
Proposition~\ref{prop:not m.5}.

We briefly explain the strategy for proving
Theorem~\ref{thm:main}(1) to clarify why this investigation will
play an important role for the proof of the main theorem. To prove
the main theorem we construct a fibred genus 2 doubly slice knot
$K$ and perform genetic modification via $\eta_i$ with $\eta_i\in
\pi_1(M_K)^{(m)}$ for all $i$. The resulting knot $K'$ is doubly
$(m)$-solvable by Corollary~\ref{cor:bi-solvable}. Then we show
that with a suitable choice of $\eta_i$ and $J_i$, for any given
double $(m)$-solution $(V_1, V_2)$ for $K'$, {\em at least one} of
$V_1$ and $V_2$ is not an $(m.5)$-solution. This will show that
$K'$ is not doubly $(m.5)$-solvable.

In fact, what we investigate was studied by Cochran and Teichner
in \cite{CT} in which they create the examples of knots that are
$(m)$-solvable but not $(m.5)$-solvable. In \cite{CT}, they show
that there is a trivial link $L\equiv \{\eta_i\}_{1\le i\le n}$
which lies in $\pi_1(M_K)^{(m)} \setminus \pi_1(M_K)^{(m+1)}$ such
that $K'$ is not $(m.5)$-solvable. However, note that to prove
Theorem~\ref{thm:main}(1) we need $K'$ to be $(n)$-solvable for
all $n$. Thus we use not the whole link $L$ but its sublinks for
genetic modification to construct our examples, and we need to
find out how to choose those sublinks.
\\
\\
\indent We follow arguments in \cite{CT}. Any result in this
section can be obtained from \cite{CT}, with a little
investigation if needed.
\\
\\
\indent Throughout this section $M$ and $M'$ denote zero surgeries
on $K$ and $K'$, respectively. We assume $\eta_i$, $J_i$, and $V$
as in the first paragraph of this section. We begin by giving a
``standard" method which gives us an $(m)$-solution $W$ for $K$
from a given $(m)$-solution $V$ for $K'$. We construct a standard
cobordism $C$ between $M$ and $M'$ as follows. For each
$(0)$-solvable knot $J_i$, choose a $(0)$-solution $W_i$ such that
$\pi_1(W_i)\cong \bbz$. We form $C$ from $M\times [0,1]$ and $W_i$
by identifying $\eta_i\times D^2$ in $M\times \{1\}$ and the solid
torus $S^1\times D^2$ in $\partial W_i = (S^3\setminus N(J_i))\cup
S^1\times D^2$ in such a way that the meridian of $\eta_i$ is
glued with the longitude of $J_i$ and the longitude of $\eta_i$ is
glued with the meridian of $J_i$ for $1\le i\le n$. ($N(J_i)$ is
an open tubular neighborhood of $J_i$ in $S^3$.) One sees that
$\partial_{-}C = M$ and $\partial_{+}C = M'$. Now we define $W$ to
be the union of the cobordism $C$ and the $(m)$-solution $V$ for
$K'$ along $M'$. Then $\partial W = M$ and $W$ is an
$(m)$-solution for $K$. To see $W$ is an $(m)$-solution for $K$,
the readers are referred to \cite{CT}.

Since $M$ fibers over $S^1$ with a fiber genus 2 closed surface
$\S$, $\pi_1(M)\cong \pi_1(\S)\rtimes\bbz$ where $\pi_1(\S)\cong
\pi_1(M)^{(1)}$. Let $S$ denote $\pi_1(\S)$. The group $S$ has a
presentation $\left<x_1,x_2,x_3,x_4\left|\right.
[x_1,x_2][x_3,x_4]\right>$. Let $(a,b)$ and $(c,d)$ be orderings
of the sets $\{1,2\}$ and $\{3,4\}$ respectively. We define the
set $P^{a,c}_n$ whose elements are pairs of elements in
$S^{(n)}$($=\pi_1(M)^{(n+1)}$) for each $n$ inductively as
follows. (Therefore we define the four sets $P^{1,3}_n, P^{1,4}_n,
P^{2,3}_n$, and $P^{2,4}_n$.) Define $P^{a,c}_1 = \{([x_a,x_b],
[x_a,x_c])_{a,c}\}$. The subscript $a,c$ for the pair is used to
designate that this pair is an element of $P^{a,c}_n$ to prevent
possible confusion in the future use. Assume $P^{a,c}_n$ has been
defined. We define $P^{a,c}_{n+1}$ as follows. For each
$(y,z)_{a,c} \in P^{a,c}_n$, $P^{a,c}_{n+1}$ contains the
following 3 pairs :
$$
([y,y^{x_a}], [z,z^{x_a}])_{a,c}, ([y,z],
[z,z^{x_a}])_{a,c}, ([y,y^{x_a}], [y,z])_{a,c}
$$
where $y^x = x^{-1}yx$. Thus $P^{a,c}_{n+1}$ has $3^{n}$ pairs.

Next, we introduce the notion of \emph{algebraic solutions}. For a
group $G$, let $G_k \equiv G/G_{tf}^{(k)}$ where $G_{tf}^{(k)}$ is
\emph{the $k^{\text{th}}$ rational derived group of $G$} by Harvey
\cite{H}. The following definition and propositions can be found
in \cite{CT}.
\begin{definition}
\label{def:algebraic solution}\cite[Definition 6.1]{CT} A
homomorphism $r : S \rightarrow G$ is called an \emph{algebraic
$(n)$-solution} ($n\ge 1$) if the following hold :
\begin{enumerate}
\item $r_\ast : H_1(S;\bbq) \rightarrow H_1(G;\bbq)$ has
2-dimensional image and there exists an ordering $(a,b)$ of the
set $\{1,2\}$ and an ordering $(c,d)$ of the set $\{3,4\}$ such
that $r_\ast(x_a)$ and $r_\ast(x_c)$ are nontrivial.
\item For each $0\le k\le n-1$, the following composition is
nontrivial even after tensoring with the quotient field
$\scrk(G_k)$ of $\bbz G_k$:
$$
H_1(S;\bbz G_k) \xrightarrow{r_\ast} H_1(G;\bbz G_k) \cong
G^{(k)}_{tf}/[G^{(k)}_{tf},G^{(k)}_{tf}] \rightarrow
G^{(k)}_{tf}/G^{(k+1)}_{tf}.
$$
\end{enumerate}
\end{definition}
We remark that if $r : S \rightarrow G$ is an algebraic
$(n)$-solution, then for any $k < n$ it is an algebraic
$(k)$-solution. The following proposition is (implicitly) proved
in the proof of Lemma $6.7$ in \cite{CT}.
\begin{prop}
\label{prop:special pair}\cite{CT} For any algebraic
$(n)$-solution $r : S \rightarrow G$ such that $r_\ast(x_a)$ and
$r_\ast(x_c)$ are nontrivial, there exists a pair in $P^{a,c}_n$
(which is called a {\bf special pair}) which maps to a $\bbz
G_n$-linearly independent set under the composition:
$$
S^{(n)} \rightarrow S^{(n)}/S^{(n+1)} \cong H_1(S;\bbz S_n)
\xrightarrow{r_\ast} H_1(S;\bbz G_n).
$$
\end{prop}

Let $W$ be the $(m)$-solution for $K$ obtained from an
$(m)$-solution $V$ for $K'$ by the ``standard" method explained as
above in this section. Let $\g \equiv \pi_1(W)^{(1)}$. The
inclusion $i : M\rightarrow W$ induces a homomorphism $h : S
\rightarrow \g$.
\begin{prop}
\label{prop:algebraic solution}\cite[Proposition 6.2]{CT} The
homomorphism $h : S \rightarrow \g$ is an algebraic
$(m)$-solution.
\end{prop}

By Proposition~\ref{prop:special pair} and
Proposition~\ref{prop:algebraic solution}, there exists an
ordering $(a,b)$ of the set $\{1,2\}$ and an ordering $(c,d)$ of
the set $\{3,4\}$ such that $h_\ast(x_a)$ and $h_\ast(x_c)$ are
nontrivial. Now we have the following proposition. We remind the
reader that $J_i$ are $(0)$-solvable and $\{\eta_i\}_{1\le i\le
n}$ is a trivial link which misses $K$. In the following
proposition, $\rho_\bbz(J_i)$ denotes the von Neumann
$\rho$-invariant $\rho(M_{J_i},\phi)$ where $\phi : \pi_1(M_{J_i})
\rightarrow \bbz$ is the abelianization. It is known that for $M$,
there is an upper bound for von Neumann $\rho$-invariants. More
precisely, there exists a constant $c_M$ such that $|\rho(M,\phi)|
\le c_M$ for every representation $\phi : \pi_1(M) \rightarrow \G$
where $\G$ is a group. (See \cite{CG} and \cite[Theorem
3.1.1]{R}.) For von Neumann $\rho$-invariants, refer to
\cite{CG,COT1,COT2}.

\begin{prop}
\label{prop:not m.5} Suppose $\rho_\bbz(J_i) > c_M$ for $1\le i\le
n$. Suppose $(a,b)$ and $(c,d)$ are orderings of the sets
$\{1,2\}$ and $\{3,4\}$ respectively such that $h_\ast(x_a)$ and
$h_\ast(x_c)$ are nontrivial in $\g$ $(= \pi_1(W)^{(1)})$. If
$\{\eta_i\}_{1\le i\le n}$ is a link in $S^3 \setminus \Sigma$
such that the set of all homotopy classes represented by $\eta_i$
contains all homotopy classes in the pairs in $P^{a,c}_{m-1}$,
then the $(m)$-solution $V$ for $K'$ is not an $(m.5)$-solution
for $K'$.
\end{prop}
\begin{proof}
By Proposition~\ref{prop:algebraic solution}, the homomorphism $h$
is an algebraic $(m)$-solution, hence an algebraic
$(m-1)$-solution. By Proposition~\ref{prop:special pair}, for the
homomorphism $h$ there exists a special pair in $P^{a,c}_{m-1}$
which maps to a $\bbz \g_{m-1}$-linearly independent set under the
composition
$$
S^{(m-1)} \rightarrow S^{(m-1)}/S^{(m)} \cong
H_1(S;\bbz[S/S^{(m-1)}]) \xrightarrow{h_\ast} H_1(S;\bbz \g_{m-1})
$$
where $\g_{m-1} = \g/\g^{(m-1)}_{tf}$. So there is at least one
pair, say $(y,z)$, in $P^{a,c}_{m-1}$ which maps to a basis of
$H_1(S;\scrk(\g_{m-1}))$ where $\scrk(\g_{m-1})$ is the (skew)
quotient field of $\bbz \g_{m-1}$. By part (2) of
Definition~\ref{def:algebraic solution}, at least one of $y$ and
$z$ maps nontrivially under the composition
$$
S^{(m-1)} \rightarrow H_1(S;\bbz \g_{m-1}) \xrightarrow{h_\ast}
H_1(H;\bbz \g_{m-1}) \rightarrow
\g^{(m-1)}_{tf}/[\g^{(m-1)}_{tf},\g^{(m-1)}_{tf}] \rightarrow
\g^{(m-1)}_{tf}/\g^{(m)}_{tf}.
$$
By our choice of $\eta_i$, this tells us that there exists
$\eta_j$ for some $j$ which maps nontrivially to
$\g^{(m-1)}_{tf}/\g^{(m)}_{tf}$, hence $i_\ast(\eta_j) \notin
\g^{(m)}_{tf} = \pi_1(W)^{(m+1)}_{tf}$.
(($\pi_1(W)^{(1)})^{(m)}_{tf} = \pi_1(W)^{(m+1)}_{tf}$ since
$H_1(W) \cong \pi_1(W)/[\pi_1(W),\pi_1(W)] \cong \bbz$ which is
torsion free.)

Let $\G \equiv \pi_1(W)/\pi_1(W)^{(m+1)}_{tf}$. Then $\G$ is an
$(m)$-solvable poly-torsion-free-abelian group by \cite[Corollary
3.6]{H}. Let $\psi : \pi_1(W) \rightarrow \G$ be the projection.
By \cite[Lemma 4.5]{CT},
$$
\rho(M,\psi|_{\pi_1(M)}) - \rho(M',\psi|_{\pi_1(M')}) =
\sum^n_{i=1}\epsilon_i\rho_\bbz(J_i)
$$
where $\epsilon_i = 0$ if $\psi(\eta_i) = e$, and $\epsilon_i = 1$
otherwise.

If $V$ were an $(m.5)$-solution for $K'$,
$\rho(M',\psi|_{\pi_1(M')}) = 0$ by \cite[Theorem 4.2]{COT1}.
Since $\psi(\eta_j) \neq e$, it follows that
$\rho(M,\psi|_{\pi_1(M)}) > c_M$, which is a contradiction.
Therefore $V$ is not an $(m.5)$-solution for $K'$.
\end{proof}

\section{The proof of the main theorem}
\label{sec:proof} We use the same notations as in
Section~\ref{sec:not m.5}. In particular, $M \equiv M_K$ and $M'
\equiv M_{K'}$. Before giving the proof, we start with our choice
for the seed knot $K$ and a little lemma for $M$. Let $T$ be the
right-handed trefoil. We define $K$ to be $T\#(-T)$. See
Figure~\ref{fig:K} below. The rectangles containing integers
symbolize full twists. Thus the rectangle labelled $+1$ symbolizes
1 right-handed full twist. Then $K$ is doubly slice by the
following theorem and its corollary due to Zeeman and Sumners,
respectively.
\begin{thm}\cite[Corollary 2, p.~487]{Z}
\label{thm:twist} Every 1-twist-spun knot is unknotted.
\end{thm}
\begin{cor}
\label{cor:twist}\cite{Su} $J\#(-J)$ is doubly slice for every
knot $J$.
\end{cor}
\noindent More generally, in \cite{Z} Zeeman proves that the
complement of a $k$-twist-spun knot in $S^4$ fibers with fiber the
punctured $k$-fold cyclic cover of $S^3$ branched along the knot
we are spinning. Also it is well-known that $J\#(-J)$ is a ribbon
knot for every knot $J$. (For instance, see \cite[Proposition 5.10
p.83]{K}.) Moreover since $T$ is a genus 1 fibred  knot, $K$ is a
genus 2 fibred knot. Combining all these, one sees that $K$ is a
genus 2 fibred doubly slice ribbon knot.

\begin{figure}[htb]
\includegraphics[scale=1.0]{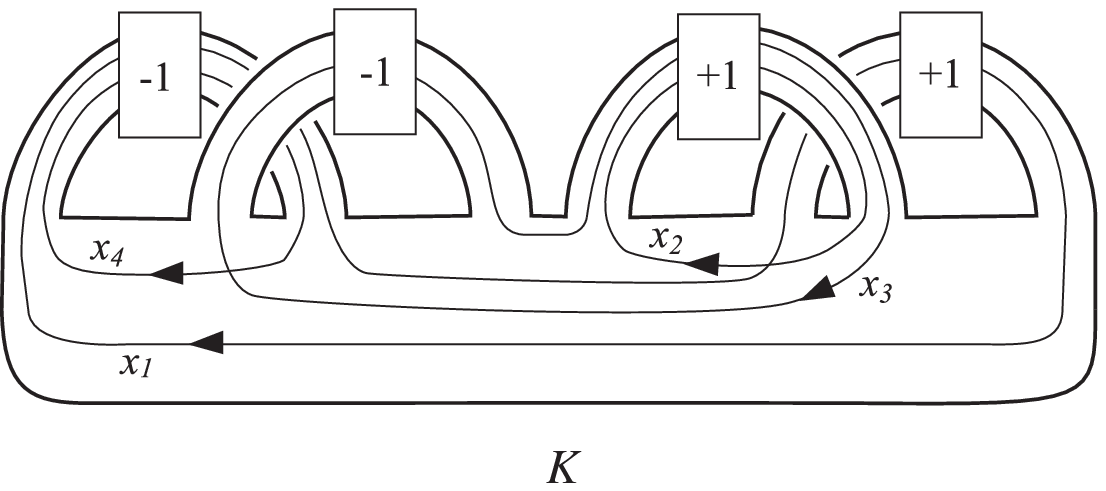}
\caption{}\label{fig:K}
\end{figure}

The knot $K$ bounds the obvious Seifert surface $F$ that is the
boundary connected sum of disks with bands as one sees in
Figure~\ref{fig:K}. Since $K$ is fibred, $M$ fibers over $S^1$
with a fiber $\S$ which is obtained by taking the union of $F$ and
a 2-disk (surgery disk) along the boundary. Let $x_1, x_2, x_3$,
and $x_4$ denote the simple closed curves on $\S$ as shown in
Figure~\ref{fig:K} whose homology classes form a symplectic basis
for $H_1(\S)$. Recall that $S = \pi_1(\S) \cong \pi_1(M)^{(1)}$.
Thus the group $S$ has a presentation
$\left<x_1,x_2,x_3,x_4\left|\right. [x_1,x_2][x_3,x_4]\right>$ as
in Section~\ref{sec:not m.5} where we abuse notations for
convenience so that each $x_i$ in the presentation is identified
with the homotopy class represented by the simple closed curve
$x_i$ on $\S$. Recall that $\L \equiv \bbz[t,t^{-1}]$.
\begin{lem}
\label{lem:generating} Any pair of $x_i$'s except for the pair
($x_1$, $x_3$) generates $H_1(M;\L)$.
\end{lem}
\begin{proof}
Denote by $y_1$ the simple closed curve which traverses once
clockwise the leftmost band on $\S$ in Figure~\ref{fig:K}.
Similarly, denote by $y_2, y_3, y_4$ the simple closed curves
traversing once clockwise the remaining bands on $\S$,
respectively. (We number $y_i$ from left to right.) Then, in
$H_1(M;\L)$, with suitable orientations for $x_i$ and $y_i$, we
have relations $x_1 = y_1 +y_4$, $x_2 = y_3$, $x_3 = y_2 + y_3$,
and $x_4 = y_1$. With the choice of basis $\{y_1,y_2,y_3,y_4\}$,
the Seifert matrix of $K$ is
\[
A = \left(\begin{matrix} -1 & 0 & 0 & 0\cr
 -1 & -1 & 0 & 0\cr
 0 & 0 & 1 & 1 \cr
 0 & 0 & 0 & 1
\end{matrix}\right).
\]
Then $H_1(M;\L)$ is presented by the matrix $tA^t - A$ with
respect to the basis $\{y_1^*, y_2^*, y_3^*, y_4^*\}$ where $A^t$
denotes the transpose of $A$ and $y_i^*$ denotes an Alexander dual
of $y_i$ in $S^3\setminus \S$. Since $A$ is invertible, $t -
A(A^t)^{-1}$ is a presentation matrix of $H_1(M;\L)$ with respect
to the basis $\{y_1, y_2, y_3, y_4\}$. Thus $H_1(M;\L)\cong
\L/(t^2-t+1)\oplus\L/(t^2-t+1)$ where $y_2$ and $y_3$ are
identified with $(1,0)$ and $(0,1)$, respectively. Also $y_1$ and
$y_4$ are identified with $(t,0)$ and $(0,t)$. Using the relations
among $x_i$ and $y_i$ and noting that $t^2-t+1=0$, one easily
deduces the lemma.
\end{proof}

Now we give the proof of Theorem~\ref{thm:main}.
\begin{proof}[Proof of Theorem~\ref{thm:main}(1)] Let $n =
2\cdot \left|P^{1,3}_{m-1}\cup P^{1,4}_{m-1}\cup
P^{2,3}_{m-1}\right| = 2\cdot 3\cdot 3^{m-2} = 2\cdot 3^{m-1}$.
(Recall that $P^{a,c}_{m-1}$ were defined in Section~\ref{sec:not
m.5}.) Let $c_M$ be a positive number given by \cite{CG} and
\cite[Theorem 3.1.1]{R} such that $\left|\rho(M,\phi)\right| \le
c_M$ for every representation $\phi : \pi_1(M) \rightarrow \G$
where $\G$ is a group. For $1\le i\le n$, let $J_i$ be an Arf
invariant zero knot such that $\rho_\bbz(J_i) > c_M$. (For
example, one can choose $J_i$ to be the connected sum of suitably
many even number of left-handed trefoils.) Since $S = \pi_1(\S)
\cong \pi_1(M\setminus \S) \cong \pi_1(S^3\setminus F)$, we can
choose $n$ simple closed curves in $S^3\setminus F$ which
represent all of the homotopy classes in the pairs in
$P^{1,3}_{m-1}\cup P^{1,4}_{m-1}\cup P^{2,3}_{m-1}$. Label these
simple closed curves by $\eta_i'$, $1\le i\le n$.

Recall that $K$ is a ribbon knot. We claim that there is a slice
disk $D$ for $K$ obtained by deforming a ribbon disk for $K$ into
$B^4$ such that $\eta_i'$ are homotopically trivial in $B^4
\setminus D$. For the proof of this claim and later use, we give
two slice disk and 4-ball pairs $(B^4, D_1)$ and $(B^4,D_2)$ for
$K$ (not $K'$) such that their union along the boundary gives an
unknotted $S^2$ in $S^4$ : from \cite{Z} and
Corollary~\ref{cor:twist}, $(B^4, D_1)$ is obtained by
half-spinning $T$ without twist and $(B^4,D_2)$ is obtained by
half-spinning $-T$ with a 1-twist. Let $W_1$ be the exterior of
$D_1$ in $B^4$ and $W_2$ the exterior of $D_2$ in $B^4$. We show
that $\eta_i'$ represent the trivial element in $\pi_1(W_1)$. This
will show the claim since $D_1$ can be obtained by deforming a
ribbon disk for $K$ into $B^4$. (To see this, refer to
\cite[Proposition 5.10 p.83]{K}.) Let $(a,b)$ be any of the
ordered pairs $(1,3)$, $(1,4)$, and $(2,3)$. Observe that the
simple closed curves $x_1$ and $x_3$ bound embedded disks in
$W_1$. (These disks are easily obtained by half-spinning without
twist the half of $x_1$ and the half of $x_2$ in $B^4$.) Therefore
$x_1 = x_3 = e$ in $\pi_1(W_1)$, hence $[x_a,x_b] = [x_a,x_c] = e$
in $\pi_1(W_1)$. Suppose $(y,z)_{a,c}$ be an element of
$P^{a,c}_j$ $(1 \le j \le m-2)$ such that $y = z = e$ in
$\pi_1(W_1)$. Then $[y,y^{x_a}] = [z,z^{x_a}] = [y,z] = e$ in
$\pi_1(W_1)$. Now using an induction argument, one sees that every
homotopy class in the pairs in $P^{a,c}_{m-1}$ represents the
trivial element in $\pi_1(W_1)$, hence $\eta_i' = e$ in
$\pi_1(W_1)$ for all $i$.

Now by Proposition~\ref{prop:ribbon} there is a trivial link
$\eta_i$ ($1 \le i \le n)$ such that each $\eta_i$ is homotopic to
$\eta_i'$ in $S^3 \setminus K$ and $K' \equiv K(\{J_i\}_{1\le i\le
n}, \{\eta_i\}_{1\le i\le n})$ is a ribbon knot. In particular,
$\eta_i$ represent all of the homotopy classes in the pairs in
$P^{1,3}_{m-1}\cup P^{1,4}_{m-1}\cup P^{2,3}_{m-1}$. Observe that
a homotopy in $S^3 \setminus K$ between $\eta_i'$ and $\eta_i$ can
be constructed by using crossing change in $S^3\setminus F$ and
the isotopy (which can be extended to the ambient isotopy). Hence
we may assume that $\eta_i$ are disjoint from $F$.

We show that $K'$ satisfies the other required conditions. To see
that $K'$ is doubly $(m)$-solvable, just observe that $\eta_i$ lie
in $\pi_1(M)^{(m)}$ which is mapped into $\pi_1(W_1)^{(m)}$ and
$\pi_1(W_2)^{(m)}$. Now it follows from
Proposition~\ref{prop:bi-solvable} that $K'$ is doubly
$(m)$-solvable.

Assume that $(V'_1,V'_2)$ is a double $(m)$-solution for $K'$. We
show that {\em at least one} of $V'_1$ and $V'_2$ is not an
$(m.5)$-solution for $K'$. Since $m\ge 2$, $(V'_1, V'_2)$ is a
double $(1)$-solution for $K'$. By Proposition~\ref{prop:1.0} and
its proof, we have $H_1(M';\L')\cong H_1(V'_1;\L')\oplus
H_1(V'_2;\L')$ where $\L' = \bbq[t,t^{-1}]$. For $i=1,2$, let
$V_i$ be the $(m)$-solution for $K$ obtained from the cobordism
$C$ and $V'_i$ as in Section~\ref{sec:not m.5}. Using the
Mayer-Vietoris sequence, one verifies that $H_1(M';\L')\cong
H_1(M;\L')$ and $H_1(V'_i;\L')\cong H_1(V_i;\L')$ for $i=1,2$. So
the inclusions $i_1 : M \rightarrow V_1$ and $i_2 : M \rightarrow
V_2$ induce the isomorphism $H_1(M;\L')\cong H_1(V_1;\L')\oplus
H_1(V_2;\L')$. We will take care of three cases : in
$H_1(V_1;\L')$, $(1)$ $(i_1)_\ast(x_1) \ne 0$, $(2)$
$(i_1)_\ast(x_3) \ne 0$, and $(3)$ $(i_1)_\ast(x_1) =
(i_1)_\ast(x_3) = 0$.

Case $(1)$ : Suppose $(i_1)_\ast(x_1) \neq 0$ in $H_1(V_1;\L')$.
Since $(i_1)_\ast$ is not a zero homomorphism (see \cite[Theorem
4.4]{COT1}), by Lemma~\ref{lem:generating} $(i_1)_\ast(x_3) \ne 0$
or $(i_1)_\ast(x_4) \ne 0$. Suppose $(i_1)_\ast(x_3) \ne 0$. Note
the homotopy classes in the pairs in $P^{1,3}_{m-1}$ are
represented by some of $\eta_i$. Thus Proposition~\ref{prop:not
m.5} implies that $V'_1$ is not an $(m.5)$-solution for $K'$. In
case $(i_1)_\ast(x_4) \ne 0$, one proves $V'_1$ is not an
$(m.5)$-solution for $K'$ using $P^{1,4}_{m-1}$ with a similar
argument.

Case $(2)$ : If $(i_1)_\ast(x_3) \ne 0$, again $V'_1$ is not an
$(m.5)$-solution for $K'$ by a reason similar to Case $(1)$. One
should use $x_1$ and $x_2$ instead of $x_3$ and $x_4$ noting the
homotopy classes in the pairs in $P^{1,3}_{m-1}$ and
$P^{2,3}_{m-1}$ are represented by $\eta_i$.

Case $(3)$ : Suppose $(i_1)_\ast(x_1) = (i_1)_\ast(x_3) = 0$. Note
$x_1 \ne 0$ and $x_3 \ne 0$ in $H_1(M;\L')$. Then $(i_2)_\ast(x_1)
\ne 0$ and $(i_2)_\ast(x_3) \neq 0$ in $H_1(V_2;\L')$ since
$H_1(M;\L')\cong H_1(V_1;\L')\oplus H_1(V_2;\L')$. By
Proposition~\ref{prop:not m.5}, since the homotopy classes in the
pairs in $P^{1,3}_{m-1}$ are represented by some of $\eta_i$,
$V'_2$ is not an $(m.5)$-solution for $K'$.

It remains to show that $K'$ is algebraically doubly slice. Using
the basis $\{x_1, x_3, x_2-x_1, x_4-x_3\}$ of $H_1(F)$, one easily
sees that the associated Seifert form of $K$ is hyperbolic. Since
$\eta_i$ are disjoint from $F$, this hyperbolic Seifert form does
not change under the above genetic modification. Hence $K'$ is
algebraically doubly slice.
\end{proof}

\begin{proof}[Proof of Theorem~\ref{thm:main}(2)]
Since a knot is $(k,\ell)$-solvable if and only if it is
$(\ell,k)$-solvable, without loss of generality we may assume
$\ell \ge k \ge 2$. Let $n = 2\cdot \left|P^{1,3}_{k-1}\cup
P^{1,4}_{k-1}\cup P^{2,3}_{k-1}\cup P^{2,4}_{\ell-1}\right| =
2\cdot (3\cdot 3^{k-2} + 3^{\ell-2})$. Let $c_M$ be the constant
as in the proof of Theorem~\ref{thm:main}(1), that is, such that
$\left|\rho(M,\phi)\right| \le c_M$ for every representation $\phi
: \pi_1(M) \rightarrow \G$ where $\G$ is a group. For $1\le i\le
n$, let $J_i$ be an Arf invariant zero knot such that
$\rho_\bbz(J_i)
> c_M$. Let $\{\eta_i\}_{1\le i\le n}$ be a trivial link in $S^3
\setminus \Sigma$ which represents all homotopy classes in the
pairs in $P^{1,3}_{k-1}\cup P^{1,4}_{k-1}\cup P^{2,3}_{k-1}\cup
P^{2,4}_{\ell-1}$. Using genetic modification, construct $K'
\equiv K(\{J_i\}_{1\le i\le n}, \{\eta_i\}_{1\le i\le n})$.

One sees that $K'$ is algebraically doubly slice using the same
reasoning as in the proof of Theorem~\ref{thm:main}(1). Let $W_1$
and $W_2$ be the slice disk exteriors for $K$ as in the proof of
Theorem~\ref{thm:main}(1). Then $K$ is $(k,\ell)$-solvable via
$(W_1,W_2)$. Since $x_1$ and $x_3$ map to the trivial element in
$\pi_1(W_1)$, the elements in the pairs in $P^{1,3}_{k-1}\cup
P^{1,4}_{k-1}\cup P^{2,3}_{k-1}$ are the trivial element in
$\pi_1(W_1)$ and in particular in $\pi_1(W_1)^{(\ell)}$. Since
$\pi_1(M)^{(\ell)}$ maps into $\pi_1(W_1)^{(\ell)}$, the elements
of $P^{2,4}_{\ell-1}$ also lie in $\pi_1(W_1)^{(\ell)}$. Regarding
$W_2$, since $\ell\ge k$, the elements of $P^{1,3}_{k-1}\cup
P^{1,4}_{k-1}\cup P^{2,3}_{k-1}\cup P^{2,4}_{\ell-1}$ lie in
$\pi_1(W_2)^{(k)}$. By Proposition~\ref{prop:bi-solvable}, $K'$ is
$(k,\ell)$-solvable (via $(W_2',W_1')$ following the notation in
Section~\ref{sec:genetic}).

Suppose $(V'_1, V'_2)$ is a $(k,\ell)$-solution for $K'$. We show
that $V'_2$ is not an $(\ell.5)$-solution for $K'$. Let $V_1$ be
the $(k)$-solution for $K$ obtained from $V'_1$ and the cobordism
$C$ between $M$ and $M'$ as in Section~\ref{sec:not m.5}. Let
$V_2$ be the $(\ell)$-solution for $K$ obtained from $V'_2$ and
$C$. As in the proof of Theorem~\ref{thm:main}(1), the inclusions
$i_1:M\rightarrow V_1$ and $i_2:M\rightarrow V_2$ induce the
isomorphism $((i_1)_\ast, (i_2)_\ast) : H_1(M;\L')
\xrightarrow{\cong} H_1(V_1;\L')\oplus H_1(V_2;\L')$. We consider
the case $k=\ell$ first. Then $P^{2,4}_{\ell-1} = P^{2,4}_{k-1}$.
Since $V_2$ is a $(k)$-solution for $K$, the inclusion $i_2 : M
\rightarrow V_2$ induces an algebraic $(k)$-solution $r_2 :
\pi_1(M)^{(1)} \rightarrow \pi_1(V_2)^{(1)}$ by
Proposition~\ref{prop:algebraic solution}. So there are orderings
$(a,b)$ and $(c,d)$ of the sets $\{1,2\}$ and $\{3,4\}$ such that
$(i_2)_\ast(x_a) \ne 0$ and $(i_2)_\ast(x_c) \ne 0$. Since all
homotopy classes in the pairs in $P^{a,c}_{k-1}$ are represented
by $\eta_i$, $V'_2$ is not a $(k.5)$- solution (i.e., not an
$(\ell.5)$-solution) for $K'$ by Proposition~\ref{prop:not m.5}.

We assume $\ell > k$. Since $V'_2$ is an $(\ell)$-solution for
$K'$, it is a $(k+1)$-solution for $K'$, so $V_2$ is a
$(k+1)$-solution for $K$. Thus $V_2$ is an algebraic
$(k+1)$-solution by Proposition~\ref{prop:algebraic solution},
hence there are orderings $(a,b)$ and $(c,d)$ of $\{1,2\}$ and
$\{3,4\}$ respectively such that $(i_2)_\ast(x_a) \ne 0$ and
$(i_2)_\ast(x_c) \ne 0$. For these $a$, $b$, $c$, and $d$, if
$(a,c)$ is one of $(1,3)$, $(1,4)$, and $(2,3)$, then since all
homotopy classes in the pairs in $P^{1,3}_{k-1}\cup
P^{1,4}_{k-1}\cup P^{2,3}_{k-1}$ are represented by $\eta_i$, by
Proposition~\ref{prop:not m.5} $V'_2$ is not a $(k+1)$-solution
for $K$, which is a contradiction. So we deduce that
$(i_2)_\ast(x_2) \ne 0$, $(i_2)_\ast(x_4) \ne 0$, and
$(i_2)_\ast(x_1) = (i_2)_\ast(x_3) = 0$. Since all homotopy
classes in the pairs in $P^{2,4}_{\ell-1}$ are represented by
$\eta_i$, by Proposition~\ref{prop:not m.5} $V'_2$ is not an
$(\ell.5)$-solution for $K'$.

Finally we show that $V'_1$ is not a $(k.5)$-solution for $K'$. If
$k=\ell$, $V'_1$ is not a $(k.5)$-solution for $K$ with the same
reason that $V_2'$ was not a $(k.5)$-solution (when $k=\ell$). If
$\ell > k$, as we showed in the previous paragraph,
$(i_2)_\ast(x_1) = (i_2)_\ast(x_3) = 0$. Since we have the
isomorphism
$$
((i_1)_\ast, (i_2)_\ast) : H_1(M;\L') \xrightarrow{\cong}
H_1(V_1;\L')\oplus H_1(V_2;\L'),
$$
it implies that $(i_1)_\ast(x_1) \ne 0$ and $(i_1)_\ast(x_3) \ne
0$. Since all homotopy classes in the pairs in $P^{1,3}_{k-1}$ are
represented by $\eta_i$, Proposition~\ref{prop:not m.5} tells us
that $V'_1$ is not a $(k.5)$-solution for $K'$.
\end{proof}

\section{bi-filtration of the double concordance group}
\label{sec:bi-filtration} We denote the double concordance group
by $\scrdc$ and the double concordance class of $K$ by $[K]$.
Since connected sum is an abelian operation, $\scrdc$ is an
abelian group. $[-K]$ is the inverse of $[K]$ in $\scrdc$ by
Corollary~\ref{cor:twist}. Recall that $K_1$ and $K_2$ are
concordant if $K_1\#(-K_2)$ is slice. Similarly, it is known that
if $K_1\#(-K_2)$ is doubly slice then $K_1$ and $K_2$ are doubly
concordant. But little is known about the double concordance group
because we have the following unanswered conjecture.

\begin{conj}
\label{conj:open conj}If knots $J$ and $K\#J$ are doubly slice,
then $K$ is doubly slice.
\end{conj}

\noindent In this section, we construct a bi-filtration of the
double concordance group using the notion of bi-solvability.

\begin{definition}
\label{def:bi-filtration} For $m,n\ge 0$, $\overline{\scrf}_{m,n}$
is defined to be the set of the double concordance classes
represented by $(m,n)$-solvable knots.
\end{definition}

\begin{prop}
\label{prop:group} $\overline{\scrf}_{m,n}$ is a subgroup of
$\scrdc$.
\end{prop}
\begin{proof}
We show that $\overline{\scrf}_{m,n}$ is closed under addition.
Let $[K_1]$ and $[K_2]$ be in $\overline{\scrf}_{m,n}$. Then
$K_1\#J_1 = K_1'\#J_1'$ and $K_2\#J_2 = K_2'\#J_2'$ for some
doubly slice knots $J_1, J_2, J_1', J_2'$ and $(m,n)$-solvable
knots $K_1',K_2'$. Thus we get $(K_1\#K_2)\#(J_1\#J_2) =
(K_1'\#K_2')\#(J_1'\#J_2')$. By Proposition~\ref{prop:connected
sum} $K_1'\#K_2'$ is $(m,n)$-solvable. Since the connected sum of
doubly slice knots is doubly slice, it follows that $K_1\#K_2$ is
doubly concordant to an $(m,n)$-solvable knot, hence $[K_1] +
[K_2] = [K_1\#K_2] \in \overline{\scrf}_{m,n}$.

Let $[K]\in \overline{\scrf}_{m,n}$. Then $K$ is doubly concordant
to some $(m,n)$-solvable knot $J$. Since $-K$ is doubly concordant
to $-J$ and $-J$ is $(m,n)$-solvable, $[-K] \in
\overline{\scrf}_{m,n}$. So the inverse of $K$ is in
$\overline{\scrf}_{m,n}$ since $-[K] = [-K]$.
\end{proof}

\begin{cor}
\label{cor:bi-filtration} $\{\overline{\scrf}_{m,n}\}_{m,n\ge 0}$
is a bi-filtration of $\scrdc$.
\end{cor}

Unfortunately, in spite of Theorem~\ref{thm:main}, it is not known
if the bi-filtration of $\scrdc$ is nontrivial because we have a
difficulty similar to Conjecture~\ref{conj:open conj}. More
precisely, it is unknown if the following is true : If $J$ and
$K\#J$ are $(m,n)$-solvable, then $K$ is $(m,n)$-solvable.

\section{doubly $(1)$-solvable knots that are not doubly
$(1.5)$-solvable}

\label{sec:1.0} 
\begin{figure}[htb]
\includegraphics[scale=1]{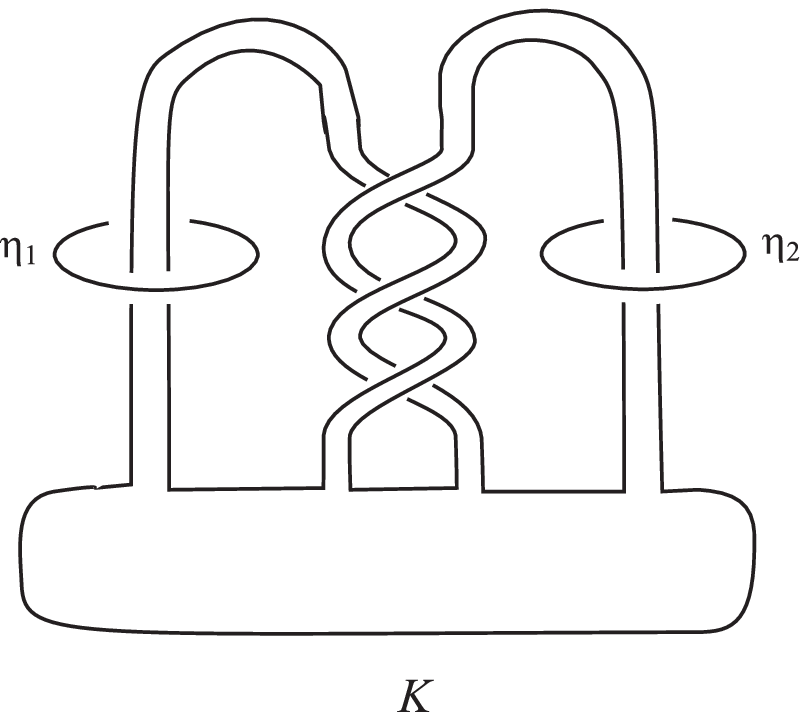}
\caption{}\label{fig:GL}
\end{figure}

In \cite{GL}, Gilmer and Livingston give a slice
knot that is algebraically doubly slice but not doubly slice.
Their example is obtained from the knot $K$ in Figure~\ref{fig:GL}
by tying the right band into a left-handed trefoil with 0-framing.
In fact, by investigating the double branched cyclic covers of
knots, they obtained an obstruction for a knot being doubly slice
in terms of the signatures of specific simple closed curves on a
Seifert surface of a knot. For more details, refer to
\cite[Theorem 4.2]{GL} and Section 5 in \cite{GL}. We prove the
following theorem.

\begin{thm}
\label{thm:1.5} There exists an algebraically doubly slice knot
$K$ that is slice and doubly $(1)$-solvable but not doubly
$(1.5)$-solvable (hence not doubly slice). Furthermore, the above
Gilmer and Livingston's obstruction vanishes for $K$.
\end{thm}

We note that $K$ in the above theorem can be shown not being
doubly slice by applying Gilmer and Livingston's method to
higher-fold finite branched cyclic covers instead of the double
branched cyclic cover.

Before proving Theorem~\ref{thm:1.0} we give useful properties of
the knot $K$ in Figure~\ref{fig:GL}. Let $a$ and $b$ be the simple
closed curves on the obvious Seifert surface $F$ which run around
the left band and the right band, respectively. With respect to
$\{a,b\}$ with a suitable choice of orientation, $K$ has the
Seifert form represented by
$$
A = \left(\begin{matrix} 0 & 1\cr 2 & 0
\end{matrix}\right)
$$
hence it is algebraically doubly slice. Since both bands are
unknotted and untwisted, $K$ is doubly slice (hence slice). Since
$tA - A^t$ is a presentation matrix of $H_1(M_K;\L')$, one sees
that $H_1(M_K;\L')\cong \L'/(t-2)\oplus \L'/(2t-1)$. That is,
there are submodules $P$ and $Q$ such that $H_1(M_K;\L') = P
\oplus Q$, and $P\cong \L'/(t-2)$, $Q\cong \L'/(2t-1)$. Here $P$
is generated by $\eta_1$ and $Q$ by $\eta_2$ where $\eta_1$ and
$\eta_2$, indicated in Figure~\ref{fig:GL}, represent the
Alexander duals of $a$ and $b$ in $H_1(S^3\setminus F)$.

Moreover the rational Blanchfield form $B\ell'$ has exactly two
self-annihilating submodules, which are $P$ and $Q$. This can be
shown easily using the presentation matrix $(1-t)(tA - A^t)^{-1}$
of $B\ell'$. Now we give the proof of Theorem~\ref{thm:1.5}. In
the proof, $\sigma_\omega$ where $\omega$ is a unit complex number
is the Levine-Tristram signature function \cite{T}. For
convenience, we define $\sigma_{r}$ $(r\in \bbq)$ to be
$\sigma_\omega$ where $\omega = e^{2\pi ir}$.

\begin{proof}[Proof of Theorem~\ref{thm:1.5}] Let $s$ be a
number such that $\frac13 < s < \frac12$. By the proof of
\cite[Theorem 1]{CL}, there exists a knot $J$ such that
$\sigma_r(J) = 0$ if $0 < r < s$ or $ 1-s < r < 1$ and
$\sigma_r(J) = 2$ if $s < r < 1-s$. Furthermore the Alexander
polynomial of $J$, say $\Delta_J(t)$, has the property that
$\Delta_J(-1) = \pm 1$ (mod 8). By \cite{leva} $J$ has vanishing
Arf invariant, and in particular $(0)$-solvable.

Let $K'\equiv K(J,\eta_2)$, the knot resulting from genetic
modification. For simplicity, let $M'\equiv M_{K'}$. Since
$\eta_2$ lies in $\pi_1(M_K)^{(1)}$, $K'$ is doubly $(1)$-solvable
by Proposition~\ref{prop:bi-solvable}. The associated Seifert form
of $F$ is hyperbolic and this Seifert form does not change under
the above genetic modification, hence $K'$ is algebraically doubly
slice. The Seifert surface $F'$ of $K'$ can be obtained from $F$
by tying the right band along $J$ with $0$-framing. Since the left
band in $F'$ remains unknotted and untwisted, $K'$ is slice. Since
$\sigma_{\frac13}(J) = \sigma_{\frac23}(J) = 0$,
$\sigma_{\frac13}(J\#J) = \sigma_{\frac23}(J\#J) = 0$. Hence the
above Gilmer and Livingston's obstruction vanishes for $K'$ (see
Theorem 4.2 and Section 5 in \cite{GL} for more details). We need
to show that $K'$ is not doubly $(1.5)$-solvable.

Since we use 0-framing when we tie a band of $F$ into $J$ to get
$F'$, $K'$ has the same Seifert matrix $A$ with respect to the
images of $a$ and $b$ under genetic modification. So
$H_1(M';\L')\cong H_1(M_K;\L')$, and $K$ and $K'$ have isomorphic
Blanchfield forms. Thus the Blanchfield form of $K'$ also has
exactly two self-annihilating submodules. For convenience we abuse
notations so that the images of $a$, $b$, $\eta_1$, $\eta_2$, $P$,
and $Q$ under genetic modification are denoted by the same
letters. So $H_1(M';\L') = P \oplus Q$, and $P\cong \L'/(t-2)$,
$Q\cong \L'/(2t-1)$. Now since $K$ is doubly slice, we have a
double $(1)$-solution $(W_1, W_2)$ for $K$ where $W_1$ and $W_2$
are the slice disk exteriors. Let $W_1'$ and $W_2'$ be the
$(1)$-solutions for $K'$ constructed as in
Section~\ref{sec:genetic}. Then $(W_1', W_2')$ is a double
$(1)$-solution for $K'$ by Proposition~\ref{prop:bi-solvable}. Let
$i_j$ be the inclusion map from $M'$ into $W'_j$ for $j = 1,2$.
Since the (rational) Blanchfield form of $K'$ has exactly two
self-annihilating submodules (which are $P$ and $Q$),
Proposition~\ref{prop:1.0} implies that we may assume $\Ker
(i_1)_\ast = Q$ and $\Ker (i_2)_\ast = P$. Since $P = P^\bot$,
$\eta_2 \notin P$, and the Blanchfield form is nonsingular, there
exists a nonzero $p \in P$ such that $B\ell(p,\eta_2) \ne 0$. By
\cite[Theorem 3.5]{COT1} $p$ induces a representation $\phi :
\pi_1(M') \rightarrow \ \Gamma_1^U$  where $\Gamma_1^U \equiv
(\bbq(t)/\L')\rtimes \bbz$. By \cite[Theorem 3.6]{COT1}, $\phi$
extends to $\psi : \pi_1(W_2') \rightarrow \Gamma^U_1$. So the von
Neumann $\rho$-invariant $\rho(M',\phi)$ can be computed using
$(W_2', \psi)$. Since $B\ell(p,\eta_2)\ne 0$, by \cite[Theorem
3.5]{COT1}  $\phi(\eta_2) \ne e$. By \cite[Proposition 3.2]{COT2}
and Property (2.2), (2.3), and (2.4) in \cite{COT2},
$$
\rho(M',\phi) = \rho(M_K,\psi|_{\pi_1(M_K)}) + \rho(M_J,
\psi|_{\pi_1(M_J)}) = \rho(M_J, \psi|_{\pi_1(M_J)}) = 2(1-2s) \ne
0.
$$

Now suppose $(V_1,V_2)$ is a double $(1)$-solution for $K'$. Let
$j_1$ be the inclusion map from $M'$ into $V_1$. Define $j_2$
similarly. Since $\Ker(j_1)_\ast$ and $\Ker(j_2)_\ast$ are
self-annihilating with respect to the rational Blanchfield form
$B\ell'$ by Proposition~\ref{prop:1.0}, without loss of generality
we may assume $\Ker (j_1)_\ast = Q$ and $\Ker (j_2)_\ast = P$. Let
$p \in P$ be as in the previous paragraph inducing the
homomorphism $\phi : \pi_1(M') \rightarrow \Gamma_1^U$. By
\cite[Proposition 3.6]{COT1} $\phi$ extends to $\psi' : \pi_1(V_2)
\rightarrow \Gamma^U_1$. So $\rho(M',\phi)$ can be computed via
$(V_2, \psi')$. If $(V_1, V_2)$ were a double $(1.5)$-solution for
$K'$, $V_2$ is a $(1.5)$-solution for $K'$. Therefore
$\rho(M',\phi) = 0$ by \cite[Theorem 4.2]{COT1}, which contradicts
the above computation that $\rho(M',\phi) \ne 0$
\end{proof}

In fact, one can show that $K'$ as above is $(1,n)$-solvable for
all $n \in \bbn$. We give another interesting example.

\begin{thm}
\label{thm:1.0} There exists a knot that is doubly $(1)$-solvable
but not $(1, 1.5)$-solvable.
\end{thm}
\begin{proof}
Let $J$ be the same as in the proof of Theorem~\ref{thm:1.5}. Let
$J_1\equiv J$ and $J_2\equiv J$. Define $K' \equiv K(\{J_1, J_2\},
\{\eta_1, \eta_2\})$, the knot resulting  from genetic
modification. Then $K'$ is doubly $(1)$-solvable but not $(1,
1.5)$-solvable. The proof follows the same course as in
Theorem~\ref{thm:1.5}, and details are omitted.
\end{proof}

{\bf Acknowledgements.} The author would like to thank Kent Orr
for his advice and helpful conversations and Peter Teichner for
essential help in showing that the main examples are ribbon knots.

\end{document}